\documentclass{amsart}
\usepackage[english]{babel}
\usepackage{graphicx}
\usepackage[hidelinks]{hyperref} 
\usepackage{amssymb}
\usepackage{amsmath}
\usepackage[foot]{amsaddr}
\usepackage{amsfonts}
\usepackage{xcolor}
\usepackage{mathtools}
\usepackage{algorithm}
\usepackage{bm}
\usepackage{booktabs} 
\usepackage{multirow,multicol}
\usepackage[noend]{algpseudocode}
\usepackage{algorithmicx, algorithm}
\usepackage[top=1in, bottom=1.25in, left=1.25in, right=1.25in]{geometry}

\usepackage{todonotes} 
\makeatletter
\algnewcommand{\LineComment}[1]{\Statex \hskip\ALG@thistlm {\color{white}\textbf{Input:}} #1}
\makeatother

\makeatletter
\algnewcommand{\LineCommentO}[1]{\Statex \hskip\ALG@thistlm {\color{white}\textbf{Output:}} #1}
\makeatother

\graphicspath{
	{Figure/heat/}
    {Figure/adv}
	{Figure/ns}
    {Figure}
}
\usepackage{float}

\begin{document}
\nocite{*}
\title[Time Extrapolation with GCA-ROM and TT decomposition]{Time Extrapolation with Graph Convolutional Autoencoder and Tensor Train Decomposition}

\author{Yuanhong Chen$^1$, Federico Pichi$^2$, Zhen Gao$^1$, Gianluigi Rozza$^2$}
\address{$^1$ School of Mathematical Sciences, Ocean University of China, Qingdao, 266100, China}
\address{$^2$ mathLab, Mathematics Area, SISSA, via Bonomea 265, I-34136 Trieste, Italy}

\begin{abstract}
Graph autoencoders have gained attention in nonlinear reduced-order modeling of parameterized partial differential equations defined on unstructured grids. Despite they provide a geometrically consistent way of treating complex domains, applying such architectures to parameterized dynamical systems for temporal prediction beyond the training data, i.e.\ the extrapolation regime, is still a challenging task due to the simultaneous need of temporal causality and generalizability in the parametric space.
In this work, we explore the integration of graph convolutional autoencoders (GCAs) with tensor train (TT) decomposition and Operator Inference (OpInf) to develop a time-consistent reduced-order model. In particular, high-fidelity snapshots are represented as a combination of parametric, spatial, and temporal cores via TT decomposition, while OpInf is used to learn the evolution of the latter. 
Moreover, we enhance the generalization performance by developing a multi-fidelity two-stages approach in the framework of Deep Operator Networks (DeepONet), treating the spatial and temporal cores as the trunk networks, and the parametric core as the branch network.
Numerical results, including heat-conduction, advection-diffusion and vortex-shedding phenomena, demonstrate great performance in effectively learning the dynamic in the extrapolation regime for complex geometries, also in comparison with state-of-the-art approaches e.g.\ MeshGraphNets.
\end{abstract}

\maketitle

\section{Introduction}
\label{sec:Introduction}
Over the last decades, reliable numerical methods, such as the Finite Element (FE), Finite Difference, and Finite Volume, have been developed to approximate the behavior of physical phenomena described by Partial Differential Equations (PDEs) \cite{quarteroni2009numerical,liu2022eighty}. 
Despite their significant advancements, when dealing with high-dimensional, parametric, time dependent and nonlinear systems, such methods encounter substantial limitations in the many-query context, since accurate prediction requires unbearable computational cost \cite{BennerSystemDataDrivenMethods2021,BennerSnapshotBasedMethodsAlgorithms2020}.

Reduced-order models (ROMs) have been developed as effective alternatives to traditional solvers to perform real-time simulations \cite{Quarteroni2016,benner2017model,hesthaven2016certified}. 
Their primary goal is to significantly enhance computational speedup while maintaining an acceptable level of accuracy, and have been widely applied across various fields, demonstrating effectiveness in diverse scenarios \cite{YoukilisProjectionbasedModelOrder2025,COPELAND2022114259,choi2021space,TenderiniModelOrderReduction2025,PichiDrivingBifurcatingParametrized2022a}. 
ROMs can be classified into intrusive \cite{S0036142900382612,17M1120531,BAI2023101969} and non-intrusive approaches \cite{hesthaven2016certified,XIAO2015522,LEE2020108973}, depending on the knowledge/exploitation of the physical equations, and implying a trade-off between interpretability, accuracy and efficiency. 
A further distinction between ROM approaches exists between linear \cite{RozzaRealTimeReduced2024,HESTHAVEN201855,GuoReducedOrderModeling2018} and nonlinear approaches \cite{LEE2020108973,TaddeiRegistrationMethodModel2020,DiezNonlinearDimensionalityReduction2021}, entailing a balance between simplicity, suitability, and generalizability.
The former methods rely on linear expansion over a set of precomputed modes, e.g.\ obtained via proper orthogonal decomposition (POD) or greedy strategies \cite{benner2017model}, while the latter often exploit autoencoder architectures, i.e.\ a nonlinear counterpart to the POD enhancing the extraction of the most meaningful features \cite{FrescaComprehensiveDeepLearningBased2021, LEE2020108973, gruberComparisonNeuralNetwork2022}.
Indeed, in the recent years, many of these strategies have been influenced by the rise of Machine Learning approaches, with a substantial shift from equation- and projection-based mathematically robust methods \cite{bergmann_enablers_2009,baiges_reduced-order_2014,akhtar_stability_2009} to data-driven and black-box models \cite{peherstorfer_dynamic_2015,brunton_data-driven_2022,FrescaPODDLROMEnhancingDeep2022}. 

Each of the aforementioned strategy has its own advantages and limitations, and the choice is usually made based on the application at hand. Well-known issues include: (i) slow Kolmogorov n-width decay jeopardizing linear compressibility \cite{ahmed_breaking_2020,noauthor_optimal_nodate}, (ii) non-affine and nonlinear problems compromising the efficiency \cite{ChaturantabutNonlinearModelReduction2010,RomorExplicableHyperreducedOrder2025,HirschNeuralEmpiricalInterpolation2024}, (iii) the extrapolation regime for time-dependent problems \cite{PeherstorferDatadrivenOperatorInference2016,lin_dynamic_2023}, (iv) lack of interpretability in black-box models \cite{RegazzoniLearningIntrinsicDynamics2024,brunton2016discovering}, (v) bifurcating scenarios with loss of uniqueness \cite{PichiArtificialNeuralNetwork2023,PichiDrivingBifurcatingParametrized2022a,PichiDeflationbasedCertifiedGreedy2025a}, (vi) ``optimal" sampling for data-driven methods \cite{MulaDynamicalApproximationSensor2025,karcher_adaptive_2022}, and (vii) complex and parameterized geometries for multi-fidelity and multi-resolution applications \cite{MORRISON2024117458,stabile_efficient_2020}.

In this manuscript, we aim at developing a data-driven strategy to comprehensively deal with parametric, nonlinear and time-dependent PDEs defined on unstructured meshes for possibly varying domains, exploiting geometric consistency and time causality while maintaining high efficiency.  

Graph neural networks (GNNs) \cite{ZHOU202057, wu2020comprehensive, Hamilton2020} have received considerable attention from the ROM community for their suitability of handling unstructured data, coming from the straightforward representation of mesh data as a graph.
Compared to traditional methods such as feedforward neural networks (FNNs) or convolutional neural networks (CNNs), GNNs demonstrate great versatility by leveraging geometric structural information, capturing complex spatial correlations, and exhibiting exceptional adaptability in processing data across unstructured meshes \cite{pfaff2020learning,han2022predicting}.

GNNs naturally emerge in the numerical solution of many PDEs and have demonstrated promising performance \cite{li_neural_2020,xiang_solving_2024,li_learning_2025,horie_physics-embedded_nodate,FrancoDeepLearningbasedSurrogate2023a,PEGOLOTTI2024107676,pfaff2020learning,gladstone2024mesh,MORRISON2024117458}. Despite being a powerful tool for learning from graph-structured data, the application of GNNs to PDE solving remains challenging due to the inherent complexity of combining varying domain geometries and physical evolution \cite{waikhom_survey_2023,ratna_inclusive_2025,liu_federated_2025,khemani_review_2024}. These challenges include: (i) computational inefficiency as GNNs are based on message-passing procedures with a rapidly increasing computational costs for large-scale problems, (ii) high memory requirements for storing and train a complete architecture, and (iii) limited interpretability as GNNs still operate as black-box models and could benefit from more meaningful physical insights.

For the dimensionality reduction of parameterized PDEs on varying geometries, the graph convolutional autoencoder (GCA) architecture has recently been proposed in \cite{PICHI2024112762}.
Falling in the class of non-intrusive and nonlinear ROM techniques, such approaches combine the versatility of graph-based architectures for unstructured data with the reduction power of autoencoders.  

A straightforward way to apply GCAs to parameterized time-dependent problems involves treating time as an additional parameter. However, this approach ignores the temporal causality inherent in such models, exhibiting a rapid accuracy degradation when performing in-time extrapolation tasks.

This issue is present for many data-driven surrogate models, which circumvent the need for explicit governing equations, but struggle to capture temporal evolution accurately.

To address this challenge and provide a more robust and generalized framework for predicting the behavior of dynamical systems, a growing body of research has focused on methods that approximate dynamical system operators from available data \cite{brunton2016discovering,PeherstorferDatadrivenOperatorInference2016}. 
Among these methods, dynamic mode decomposition (DMD) \cite{schmid2010dynamic} and its extensions, are effective in extracting linear approximation models and are widely recognized for their simplicity and ease of implementation. 
However, such linear approximation may fail to capture complex dynamic features adequately, particularly when dealing with nonlinear dynamical systems. 
While the combination of DMD with Koopman operators enables the transformation of finite-dimensional nonlinear systems into infinite-dimensional linear systems \cite{brunton2016koopman,rowley2009spectral,Williams}, selecting appropriate observation functions remains a challenge. Operator inference (OpInf) \cite{PeherstorferDatadrivenOperatorInference2016, mcquarrie2023nonintrusive}, on the other hand, constructs ROMs with polynomial nonlinearities in a non-intrusive way. Unlike traditional projection-based methods that require access to and modification of the full-order model, OpInf uses least-squares regression to directly infer low-dimensional operators from high-fidelity simulation data. The resulting ROMs can then be solved efficiently to generate approximate solutions.
Notably, in cases where the constructed ROMs system is linear and inputs are time-invariant, the operators derived through OpInf are equivalent to those obtained using DMD.

Furthermore, for high-dimensional tensors containing the information from several parameters, conventional dimensionality reduction methods such as POD often require unfolding or matricization, which disrupts the inherent structure of the tensor. In contrast, tensor decomposition approaches provides an effective strategy for obtaining low-rank approximations of high-dimensional tensors while preserving their original structure. Tensor Train (TT) decomposition  \cite{Oseledets2011} is a commonly used method that significantly reduces memory requirements and computational complexity by decomposing a \(d\)-dimensional tensor into a sequence of \(d\) low-rank core tensors.  Compared to traditional methods such as canonical polyadic decomposition \cite{Hitchcock1927} and Tucker decomposition \cite{Tucker1966}, TT decomposition is more computationally efficient and effectively mitigates the NP-hardness issue associated with tensor decomposition \cite{Li2022a}.

This study aims to leverage the strengths of GCA-ROM as a surrogate model in addressing parameter-dependent PDEs defined on complex and varying geometries, and extend its applicability to dynamical systems.

In particular, the main contribution of this paper is the development of a surrogate low-fidelity strategy, namely LF-TTOI, exploiting GCA-ROM for geometric consistency, TT decomposition for dimensionality reduction and OpInf for time extrapolation.
In the offline stage, the TT decomposition is performed on the sampled high-fidelity snapshots, collected for different parameters and times instances, to obtain parameter-, spatial-, and time-dependent cores, enabling the separate treatment of varying features, so that they can be used to reconstruct the solutions. 
Exploiting the time core obtained from the TT decomposition, an OpInf model is built to simulate the temporal evolution.
Moreover, a number of GCA-ROMs are trained at selected time points to compute the coefficients for new parameters at given time slices.

During the online phase, to approximate the solution for unseen time/parameter values, we combine the time core predicted from the OpInf model, and the parametric core obtained by minimizing the error between GCA-ROM and the low-rank approximation. 

Building on this, we also propose a multi-fidelity framework from the perspective of DeepONet \cite{Lu2019LearningNO}. 
The idea behind the developed MF-TTOI-DeepONet approach is to augment the amount of available information by incorporating the predicted coefficients from the LF-TTOI as low-fidelity approximations, and exploit the branch and trunk network structure of DeepONet to separately treat the different cores.

In this hybrid framework, the branch network benefits from such augmented input comprising the parameter and the LF-TTOI predictions for residual learning, while the trunk network is simply represented by the combination of time and spatial cores.
Moreover, to fully leverage data of different fidelities, the training of the branch network is conducted in two stages: a pre-training stage based on GCA-ROM estimates and a fine-tuning stage using high-fidelity snapshots. 
MF-TTOI-DeepONet can be seen as a further enhancement of the LF-TTOI strategy, requiring an additional training of the branch network during the offline step, but also providing an improved reconstruction accuracy by exploiting all the available information at hand. 

The proposed approaches are evaluated on benchmark problems including: (i) heat conduction, (ii) advection-diffusion, and (iii) vortex shedding. Results demonstrate their performance in accurately resolving the dynamical behavior of complex fields on parameterized geometries, achieving great and robust performance even within the scarce data context and in the extrapolation regime. Furthermore, the proposed framework has been also compared with state-of-the-art intrusive and non-intrusive methods, such as: projection-based POD (POD-Proj), POD-Neural Network \cite{HESTHAVEN201855}, and MeshGraphNets \cite{pfaff2020learning}, showing great reliability and generalization properties in complex scenarios.

The remainder of this paper is organized as follows.
Section \ref{Problem setup} presents the problem formulation to establish the theoretical foundation of the study, followed by a review of related methodologies.
Section \ref{Proposed} detailed explained the proposed methods. 
The results obtained for different test cases are discussed in Section \ref{Numerical examples}.
Finally, Section \ref{Conclusion} concludes the paper with a summary of findings and future perspectives.

\section{Computational framework and related methodologies}

In this section we aim at defining the main class of problems we are interested in, and the methodologies that constitute the backbone of our approach. In particular, we start by introducing the general formulation for parameterized, nonlinear and time-dependent PDEs defined on possibly varying domains. Then, we briefly review the architectural choices related to the parameterized geometries (GCA-ROM), the dimensionality reduction (TT decomposition), and the extrapolation in-time (OpInf).

\subsection{Problem formulation}
\label{Problem setup}
Consider the general form of a parameterized time-dependent PDE expressed as:
\begin{equation}
	\label{PDE}
	\begin{cases}
		\dfrac{\partial u(\bm{x}, t; \bm{\mu})}{\partial t} = f(\bm{x}, t,u(\bm{x}, t; \bm{\mu});\bm{\mu}), & (\bm{x}, t) \in \Omega_{\bm{\mu}} \times [0,t_f],\\
		\mathcal{B}(\bm{x}, t,u(\bm{x}, t; \bm{\mu});\bm{\mu})= b(\bm{x}, t,u(\bm{x}, t; \bm{\mu});\bm{\mu}), & (\bm{x}, t) \in \partial \Omega_{\bm{\mu}}\times  [0,t_f], \\
		u(\bm{x}, 0; \bm{\mu}) = u_0(\bm{x}; \bm{\mu}), & \bm{x} \in \Omega_{\bm{\mu}},
	\end{cases}
\end{equation}
where \( u(\bm{x}, t; \bm{\mu}) \) is the unknown state of the governing equation, \( \bm{x} \) denotes the spatial coordinates defined on a reference configuration of the parameterized computational domain \( \Omega_{\bm{\mu}} \), \( t \) is the time variable, and \( \bm{\mu} \in \mathcal{P} \) represents the parameter vector describing the features of the system. In our setting, the parameters \( \bm{\mu} \) can encompass a wide range of physical or geometric properties, such as material coefficients (e.g., conductivity, viscosity, or diffusivity), boundary conditions, source terms, or characteristic geometric attributes of the domain. The function $f(\bm{x}, t, u(\bm{x}, t; \bm{\mu});\bm{\mu})$ denotes the PDE operator containing the spatial differential operators and the source term. Finally, $u_0(\bm{x}; \bm{\mu})$ and $\mathcal{B}(\bm{x}, t, u(\bm{x}, t; \bm{\mu});\bm{\mu})$ represent the initial and boundary conditions, respectively.

After spatial discretization using a suitable high-fidelity approximation, e.g.,\ the FE method \cite{quarteroni2009numerical}, the PDE system in Equation \eqref{PDE} is transformed into a semi-discrete form involving the state vector $\mathbf{u}_{\mathcal{N}}(t;\bm{\mu}) \in \mathbb{R}^{N_h}$, expressed as:
\begin{equation}
	\label{dynamical-system}
	\begin{cases}
		\dfrac{\mathrm{d} \mathbf{u}_{\mathcal{N}}(t;\bm{\mu})}{\mathrm{d} t} = \bm{f}(t,\mathbf{u}_{\mathcal{N}}(t;\bm{\mu});\bm{\mu}), \\
		\mathbf{u}_{\mathcal{N}}(t_0; \bm{\mu}) = \mathbf{u}_{\mathcal{N}}^0(\bm{\mu}),
	\end{cases}
\end{equation}
where $N_h$ is the number of degrees of freedom after discretization. The function $\bm{f}$ represents the discrete counterpart of the PDE operator, while $\mathbf{u}_{\mathcal{N}}^0(\bm{\mu})$ corresponds to the discretized initial condition.

High-fidelity models provide accurate solutions but are often characterized by the curse of dimensionality, resulting in computationally expensive simulations.
Thus, numerically approximating such problems becomes computationally unbearable, especially when repeated evaluations are required for varying parameter configurations or when the solution needs to be retrieved in a short amount of time, i.e.\ the many-query and real-time contexts.
These challenges motivated the development of ROMs and, more in general, surrogate models aimed at providing efficient and accurate approximations for a wide class of academic and engineering applications \cite{Quarteroni2016,benner2017model,hesthaven2016certified}. We remark that, contrarily to common machine learning applications, the main starting point of ROMs is the availability of a limited computational budget, i.e.\ a scarce amount of high-fidelity data, which is used to train the surrogate model. Thus, we restrict our attention to the scarce-data regime, where only a small portion of the sampled dataset is used for training the model, while the remaining part is employed to test its generalization capabilities.  
In this context, we focus on constructing a data-driven ROM based on the following ingredients: (i) GCA-ROM to include geometric biases in the training phase for parameterized geometries, (ii) TT decomposition to decouple the contribution of the input features, and (iii) OpInf to learn a dynamical model for extrapolating in-time the dynamics.

\subsection{Graph convolutional autoencoder for nonlinear model order reduction}\label{gnn}
GNNs have emerged as a powerful and versatile framework for processing graph-structured data \cite{ZHOU202057, wu2020comprehensive, Hamilton2020}. They can leverage the features of individual nodes/edges and their structural relationships, e.g.\ connection between neighbors and sparsity, enabling significant advancements across various applications.

In the context of model order reduction, many contributions have been recently presented following the seminal work of Pfaff et al.\ \cite{pfaff2020learning} introducing the MeshGraphNet architecture.

MeshGraphNet employs a GNN architecture to predict system dynamics at successive time instances based on current and historical mesh states, using learned message-passing on graph-structured meshes. An enhanced GNN \cite{gladstone2024mesh}, building on MeshGraphNets, is developed to overcome the challenge of capturing long-range interactions in static solid mechanics.  In \cite{10913800,ALRASHDI2025143}, a Long Short-Term Memory (LSTM) network is integrated in the GNN framework to capture long-term dependencies, while the Deep Graph Operator Network framework \cite{sun2023deepgraphonet} combines graphs with DeepONet to approximate the dynamics of complex systems with underlying sub-graph structures. In \cite{PEGOLOTTI2024107676}, a GNN-based surrogate model is trained on three-dimensional hemodynamic data, accurately predicting pressure and flow rates along vessel centerlines with high computational efficiency.

Building on the strengths of GNNs, GCA-ROM has recently gained attention for its effectiveness in nonlinear model order reduction \cite{PICHI2024112762}, being particularly suitable for PDEs defined on unstructured and parameterized grids, as it employs an autoencoder architecture based on Graph Convolutional layers. 
In contrast to traditional CNNs, which operate on grid-structured data such as images, convolutional operations and pooling operators defined on graph data ensure permutation invariance, meaning the output remains consistent regardless of the ordering of nodes. 

The key operation of graph convolutions lies in aggregating information from neighboring nodes through a message-passing mechanism. 
At each message-passing step \( k \), the features of the \( \mathcal{N}(u) \) neighboring nodes of a given node \( u \) are aggregated and combined to update the node's feature representation, enabling the propagation of information across the graph structure. 

The graph convolutional encoder projects high-dimensional system states onto a low-dimensional latent space, while the decoder reconstructs the original high-dimensional data from the latent representation. This way, not only the physical features are processed, but also the geometric knowledge is directly incorporated into the model.
Finally, by learning the relationship between parameter and the latent space via an auxiliary and small neural network, GCA-ROM enables the decoder to reconstruct high-dimensional targets from the predicted latent representation corresponding to the new parameter value, directly into its reference geometrical configuration. 

\begin{figure}[htbp]
\begin{center}
    \includegraphics[width=0.8\textwidth]{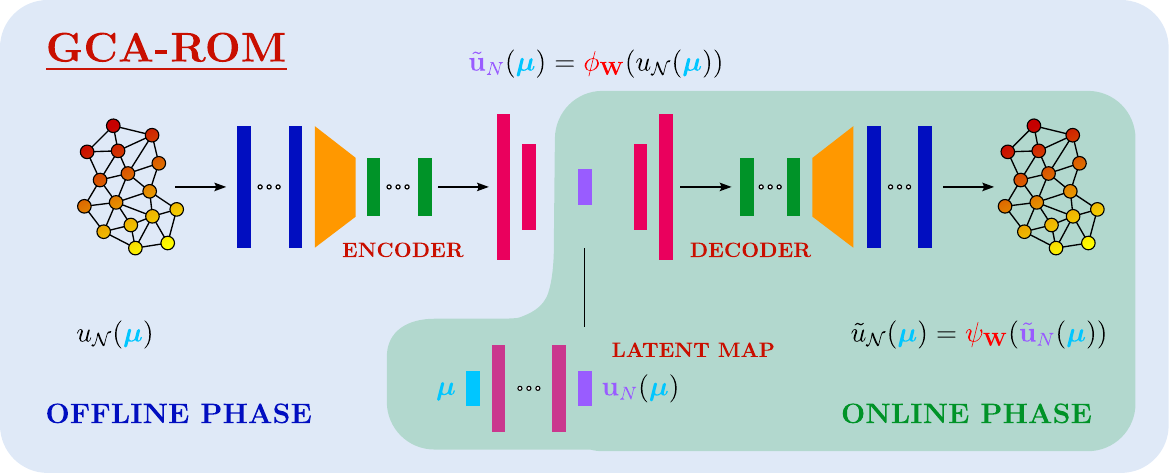}
\end{center}
\caption{Visual representation of the GCA-ROM.}
\label{fig:flow}
\end{figure}

For clarity, we report the structure of GCA-ROM in Figure \ref{fig:flow}. Specifically, we start by considering a set of training parameters, denoted as $\Xi_{\mathrm{train}} = \{\bm{\mu}_{\text{train}}^1,\bm{\mu}_{\text{train}}^2,\ldots,\bm{\mu}_{\text{train}}^{N_{\text{train}}}\}$,
and for each parameter $\bm{\mu} \in \Xi_{\mathrm{train}}$, we compute a FE high-fidelity solution $\mathbf{u}_{\mathcal{N}}(\bm{\mu})$ defined on the parameterized domain $\Omega(\bm{\mu})$.

Exploiting the graph dataset $\{\bm{\mu}_{\text{train}}^{i}, \mathbf{u}_{\mathcal{N}}(\bm{\mu}_{\text{train}}^{i}), \Omega(\bm{\mu}_{\text{train}}^{i})\}_{i=1}^{{N_{\text{train}}}}$, GCA-ROM constructs an autoencoder that approximates the identity map via two nonlinear mappings: an encoder $\tilde{\mathbf{u}}_N(\boldsymbol{\mu}) = \phi_{\mathbf{w}}(\mathbf{u}_{\mathcal{N}}(\boldsymbol{\mu}))$ and a decoder $\tilde{\mathbf{u}}_{\mathcal{N}}(\boldsymbol{\mu}) = \psi_{\mathbf{w}}(\tilde{\mathbf{u}}_N(\boldsymbol{\mu}))$, such that
\[
\tilde{\mathbf{u}}_{\mathcal{N}}(\boldsymbol{\mu}) = \psi_{\mathbf{w}}( \phi_{\mathbf{w}}(\mathbf{u}_{\mathcal{N}}(\boldsymbol{\mu}))),
\]
where $w$ denotes the weights of the architecture, and we denoted by $\tilde{\mathbf{u}}_{\mathcal{N}}(\boldsymbol{\mu})$ the reconstructed solution, using a FNN to learn a mapping $\boldsymbol{\mu} \mapsto {\mathbf{u}}_N(\boldsymbol{\mu})$ from parameters to the latent space during the offline phase. 

Different choices can be made to define the encoder/decoder modules, such as the type of graph convolutional layers, the strategy to up- and down-sample the mesh, and the residual learning approach. Despite that, an ablation study showed a nice robustness in terms of performance of the overall architecture, even in the data scarcity context.
More details concerning all the components of the modular architecture can be found in \cite{PICHI2024112762}.

To evaluate the surrogate model for any given parameter $\bm{\mu} \in \mathcal{P}$, a regression task is performed during the offline phase for the latent mapping, thereby enabling the reconstruction of solutions at arbitrary parameters without the need for encoding.  
The autoencoder and the FNN are trained simultaneously, and the input to the decoder is the output of the multilayer perceptron to enhance the precision of latent variable estimation during network optimization, i.e.\ $\tilde{\mathbf{u}}_{\mathcal{N}}(\boldsymbol{\mu}) = \psi_{\mathbf{w}}({\mathbf{u}}_N(\boldsymbol{\mu}))$.
Finally, the total loss function is defined as
\begin{equation}
\mathcal{L}=\frac{1}{N_{\mathrm{train}}} \sum_{i=1}^{N_{\mathrm{train}}}\left\|\tilde{\mathbf{u}}_{\mathcal{N}}\left(\bm{\mu}_{\text{train}}^i\right)-\mathbf{u}_{\mathcal{N}}\left(\bm{\mu}_{\text{train}}^i\right)\right\|_2^2 + \frac{\lambda}{N_{\mathrm{train}}}  \sum_{i=1}^{N_{\mathrm{train}}}\left\|\tilde{\mathbf{u}}_N\left(\bm{\mu}_{\text{train}}^i\right)-{\mathbf{u}}_N\left(\bm{\mu}_{\text{train}}^i\right)\right\|_2^2,
\end{equation}
where $\lambda$ is the weight balancing the two components. 

Although GCA-ROM is effective in approximating the solutions for PDEs characterized by physical and geometric parameters, the original architecture lacks of consistency when dealing with dynamical systems. Indeed, when applied to time-dependent problems, GCA-ROM is agnostic to the correlation between successive time-steps, with no preservation of the temporal causality inherent in the system.   
Two possible strategies are available to directly apply graph autoencoders to time-dependent problems. The first is to treat time $t$ as an additional parameter and then perform the standard GCA-ROM procedure, thus neglecting the evolution of the dynamics.
The second is to perform a state-to-state prediction, whereby the next state of the system is predicted based on the current state. 
In both cases, the training process are computationally much more expensive as they require a large amount of snapshots from each trajectory to accurately describe the behavior in time.

To address these challenges, and fully exploit the strengths of GCA-ROM, we integrate TT decomposition and OpInf to effectively capture the dynamics of the system in the extrapolation regime. 

\subsection{Tensor train decomposition}
A key ingredient for our approach is the TT decomposition \cite{Oseledets2011,10930561}, consisting of a strategy to provide compact representation of high-dimensional tensors, here resulting from the collection of spatial snapshots in time and parameter space.
In general, for a \(d\)-dimensional tensor $\mathcal{K} \in \mathbb{R}^{n_1\times n_2\times\cdots\times n_d}$, where \(n_i\) denotes the size of the $i$-th dimension, the decomposition approximates the tensor $\mathcal{K}$ as a combination of $d$ smaller second-order or third-order tensors. 
These low-order tensors are sequentially connected through contraction operations, resembling the arrangement of train cars, to obtain an approximation $\bm{K} \in \mathbb{R}^{n_1\times n_2\times\cdots\times n_d}$ given by:
\begin{equation}
	\bm{K} = \sum_{\alpha_0 = 1}^{r_0} \cdots \sum_{\alpha_d = 1}^{r_d} \bm{G}_1(\alpha_0, :, \alpha_1)\otimes \bm{G}_2(\alpha_1, :, \alpha_2)\otimes \cdots \otimes \bm{G}_d(\alpha_{d-1}, :, \alpha_d),
\end{equation}
where the tensors $\bm{G}_k \in \mathbb{R}^{r_{k-1} \times n_k \times r_k}$ for $k = 1, \dots, d$ are referred to as the cores,  the sequence $\{r_i\}_{i=0}^d$ denotes the TT ranks \cite{Oseledets2011}, with $r_0 = r_d = 1$ corresponding to the two second-order initial and final tensors, and $\otimes$ is the Kronecker product of matrices.

By processing smaller tensors, high-dimensional tensors can be efficiently represented and indirectly handled, thereby reducing the computational complexity and storage requirements. 
TT decomposition is typically achieved through \(d-1\) singular value decompositions (SVDs), and the TT ranks are determined by appropriately truncating the singular values during each SVD step. 
After performing TT decomposition, the final approximation error is bounded in the Frobenius norm as follows:  
\[
\| \mathcal{K} - \bm{K}\|_F^2 \leq \sum_{k=1}^{d-1} \varepsilon_k^2,
\]
where \(\varepsilon_k\) represents the truncation error from the \(k\)-th SVD step.

In the context of parameterized time-dependent problems, we aim at decomposing the high-dimensional tensor obtained from the solution snapshots w.r.t.\ its spatial, parametric and temporal features, i.e.\ $d=3$. Having denoted the parameter set by $\Xi_{\mathrm{train}}$, for each parameter $\bm{\mu} \in \Xi_{\mathrm{train}}$, we collect $N_t$ high-fidelity solutions $\{\mathbf{u}_{\mathcal{N}}(\bm{\mu},t^i)\}_{i=1}^{{N_t}}$ defined on the spatial domain  $\Omega(\bm{\mu})$ and store the tensor $\mathbf{S}_{\text{hf}} \in \mathbb{R}^{N_{\text{train}}\times N_h\times N_t}$. 
Performing the TT decomposition on $\mathbf{S}_{\text{hf}}$, we obtain the following approximation: 
$$
\mathbf{S}_{\text{hf}}\approx \sum_{\alpha_0 = 1}^{r_0} \cdots \sum_{\alpha_d = 1}^{r_{3}} \mathbf{G}_1(\alpha_0, :, \alpha_1)\otimes \mathbf{G}_2(\alpha_1, :, \alpha_2)\otimes \mathbf{G}_{3}(\alpha_{2}, :, \alpha_{3}),
$$
where $\mathbf{G}_1 = \left[\bm{g}_{1}^{1}, \bm{g}_{1}^{2}, \cdots, \bm{g}_{1}^{N_{\text{train}}}\right] \in \mathbb{R}^{N_{\text{train}} \times r_1}$ represents the parameter-dependent core, 
$\mathbf{G}_2 \in \mathbb{R}^{r_1 \times N_h \times r_2}$ corresponds to the space-dependent core, and $\mathbf{G}_3 = \left[\bm{g}_{3}^{1}, \bm{g}_{3}^{2}, \cdots, \bm{g}_{3}^{N_t}\right] \in \mathbb{R}^{r_2 \times N_t}$ represents the  time-dependent core. We remark that $\bm{g}_{1}^{i} \in \mathbb{R}^{1 \times r_1}$ and $\bm{g}_{3}^{i} \in \mathbb{R}^{r_2 \times 1}$ are respectively parameter- and time-dependent, and we exploit the latter to model the temporal evolution of such low-dimensional time coefficients employing OpInf.

\subsection{Operator inference}\label{sec:opinf}
Operator inference is a data-driven framework for learning operators that model the dynamics of a system, aiming to infer a polynomial model describing the temporal evolution of the state \(\boldsymbol{g}(t) \in \mathbb{R}^{r_2 \times 1}\) from data \(\{\boldsymbol{g}(t^i)\}_{i=1}^{N_t}\) collected at \(N_t\) evenly spaced time steps.

The inferred model is commonly expressed as:  
\begin{equation}
\label{op}
\frac{\mathrm{d} \boldsymbol{g}(t)}{\mathrm{d} t} = \sum_{\ell=1}^L \boldsymbol{A}_{\ell} \boldsymbol{g}^{\ell}(t) + \boldsymbol{B}\boldsymbol{u}(t)+\boldsymbol{c},
\end{equation}
where $L$ is the highest order of the polynomial nonlinearity, \(\boldsymbol{A}_{\ell} \in \mathbb{R}^{r_2 \times n_{\ell}}\), \(\boldsymbol{B} \in \mathbb{R}^{r_2 \times p}\), and $\boldsymbol{c}\in \mathbb{R}^{r_2 \times 1}$  are the inferred operators, \(\boldsymbol{g}^{\ell}(t) \in \mathbb{R}^{n_{\ell} \times 1}\) represents the \(\ell\)-th order polynomial term of the state vector with \(n_{\ell} = \binom{r_2 + \ell - 1}{\ell}\) for \(\ell = 1, \ldots, L\),  and \(\boldsymbol{u}(t)\in  \mathbb{R}^{p \times 1}\) represents the external input.

Derivative information of the state vector is crucial for inferring the underlying system dynamics. When direct derivative data is unavailable, the time derivatives of the states, \(\dot{\boldsymbol{g}}(t_j)\), at time \(t_j\) for \(j = 0, \ldots, N_t\), can be approximated using methods such as finite difference, yielding the following approximation $\boldsymbol{Q} = \left[\dot{\boldsymbol{g}}(t^1), \ldots, \dot{\boldsymbol{g}}(t^{N_t})\right]$. 

With the available data \(\{\boldsymbol{g}(t^i), \dot{\boldsymbol{g}}(t^i), \mathbf{u}(t^i)\}_{i=1}^{N_t}\), representing the states,  time derivatives, and inputs at \(N_t\) time instances, the operators \(\boldsymbol{A}_\ell\), for \(\ell = 1, \ldots, L\), and \(\boldsymbol{B}\) can be identified by solving the following least-squares optimization problem:  
\[
\min_{\boldsymbol{A}_1, \ldots, \boldsymbol{A}_L, \boldsymbol{B}} \frac{1}{N_t} \sum_{k=1}^{N_t} \left\| \sum_{\ell=0}^L \boldsymbol{A}_{\ell} \boldsymbol{g}_{k}^{\ell} + \boldsymbol{B} \boldsymbol{u}_k + \boldsymbol{c} - \dot{\boldsymbol{g}}_k \right\|_2^2,
\]  
where \(\boldsymbol{g}_{k}^{\ell}\) represents the \(\ell\)-th order polynomial term of the state vector \(\boldsymbol{g}_k = \boldsymbol{g}(t^k)\), and \(\dot{\boldsymbol{g}}_k = \dot{\boldsymbol{g}}(t^k)\) is the time derivative of the state vector at time \(t^k\). 
This can also be expressed in linear form as follows:  
\begin{equation}
\label{min}
\min_{\boldsymbol{{O}}} \left\|\boldsymbol{D} \boldsymbol{{O}} - \boldsymbol{R}\right\|_F^2,
\end{equation}
where \(\boldsymbol{D}= \left[\mathbf{1}_{N_t}, \boldsymbol{G}_1^T, \ldots, \boldsymbol{G}_L^T, \boldsymbol{U}^T\right]\in \mathbb{R}^{N_t \times \bar{n}}\), \(\boldsymbol{O} = \left[\boldsymbol{c}, \boldsymbol{A}_0, \ldots, \boldsymbol{A}_L, \boldsymbol{B}\right]^T \in \mathbb{R}^{\bar{n} \times r_2}\), and \(\boldsymbol{R} \in \mathbb{R}^{N_t \times r_2}\) record the derivative information. Here, $\mathbf{1}_{N_t}\in \mathbb{R}^{N_t}$ is a vector of ones, $\boldsymbol{G}_i = \left[\boldsymbol{g}_{1,i}, \ldots, \boldsymbol{g}_{N_t,i}\right]$,  $\boldsymbol{U} = \left[\boldsymbol{u}_1, \ldots, \boldsymbol{u}_{N_t}\right]$, and \(\bar{n} = 1 + \sum_{\ell=1}^L n_{\ell}+p\).

The minimization problem in Equation \eqref{min} can be solved using standard linear algebra techniques. However, due to the insufficiency of observations, regularization is necessary to enhance the formulation robustness and algorithmic scalability of the approach. Thus, the OpInf method aims to solve the following optimization problem:  
\begin{equation}
\label{ls}
\min_{\boldsymbol{{O}}} \left\|\boldsymbol{D} \boldsymbol{{O}} - \boldsymbol{R}\right\|_F^2 + \mathcal{R}(\boldsymbol{{O}}),
\end{equation}
where \(\mathcal{R}(\boldsymbol{{O}})\) is a regularization term that penalizes the entries of the learned operators. Once this minimization problem is solved, the evolution of the state vector can be explicitly represented, enabling the prediction of future states.

\section{Graph-based approaches for time-dependent problems}\label{Proposed}
In this section, we aim to develop data-driven surrogate models based on GCA-ROM, TT and OpInf to accurately and efficiently investigate parameterized time-dependent PDEs defined on complex geometries.  
From the above discussion, it is clear that GCA-ROM provides an effective strategy for nonlinear model order reduction by leveraging graph convolutional layers for autoencoder architectures.
However, its practical applicability is limited by the need for high-fidelity data across the entire time domain and the associated memory storage and computational costs required for training.  

To make GCA-ROM applicable to time-dependent problems, we propose a novel low-fidelity framework exploiting TT decomposition and operator inference, which consistently treats time and geometrical parameters efficiently. Moreover, inspired by the DeepONet's framework, we extend our strategy in a multi-fidelity context, to further enhance predictive accuracy by effectively leveraging both low- and high-fidelity data.

\subsection{Low-fidelity approach with TT decomposition and OpInf}
As discussed in Section~\ref{sec:opinf}, the temporal evolution of the time-dependent coefficients $\bm{g}_{3}(t)$ can be approximated as $\tilde{\bm{g}}_{3}(t)$ by applying the OpInf technique to $\{\bm{g}_{3}^{i}\}_{i=1}^{N_t}$ extracted from the core tensor $\mathbf{G}_3$.
This is particularly significant because, in many scenarios, the governing dynamical equations of the system are unknown. 

Moreover, obtaining comprehensive data across the entire temporal domain is often challenging. Thus, once we have learned the evolution of the time-dependent core, we can exploit the surrogate model to perform temporal extrapolation, i.e., for each parameter $\bm{\mu}_{\text{train}}^i \in \Xi_{\mathrm{train}}$, reconstruct the solution in the form $\boldsymbol{u}_{\mathcal{N}}(\bm{\mu}_{\text{train}}^i, t) \approx \bm{g}_{1}^{i}\otimes\bm{G}_{2}\otimes \tilde{\bm{g}}_{3}(t)$ for time  $t\in [0, t_f]$ in the case where only snapshots in $[0, \mathcal{T}]\subseteq[0, t_f]$ are available. 

\begin{figure}[htbp]
\begin{center}
    \includegraphics[width=0.9\textwidth]{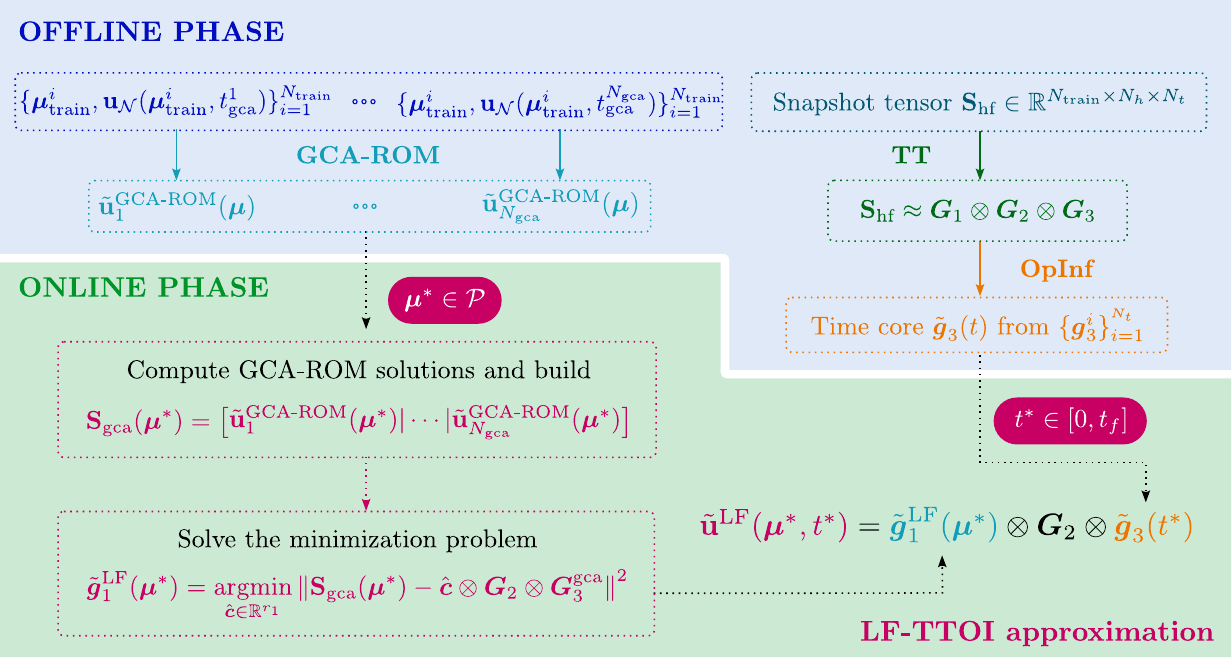}
\end{center}
\caption{The framework of LF-TTOI.}
\label{lf}
\end{figure}

However, our objective is to also apply the proposed method for the prediction of parameterized nonlinear systems, which requires computing the predicted coefficients $\tilde{\bm{g}}_{1}(\bm{\mu}^*)$ for any parameter $\bm{\mu}^* \in \mathcal{P}$. For this purpose, a fast computation method for the parametric core is needed. We propose to combine GCA-ROM with TT decomposition and OpInf, namely the LF-TTOI approach, where low-fidelity stands for the surrogate strategy exploited for the coefficients, and that will be compared to the multi-fidelity approach discussed in the next section. 
The proposed LF-TTOI procedure is detailed in Algorithm \ref{LF-TTOI}, while a schematic representation of the framework is shown in Figure \ref{lf}.

\begin{algorithm}[t]
\caption{LF-TTOI}
\label{LF-TTOI}
\vskip+10pt
\begin{algorithmic}[1]
\Statex{{\bf \underline{Offline Phase}}} 
\Statex{{\bf Input:} Training dataset with parameter samples $\Xi_{\mathrm{train}}=\{\bm{\mu}_{\text{train}}^i\}_{i=1}^{N_{\text{train}}}$, time steps $\{t^j\}_{j=1}^{N_t}$,}
\LineComment{GCA time points $\{t_{\text{gca}}^j\}_{j=1}^{N_{\text{gca}}}$, and high-fidelity solutions $\{\mathbf{u}_{\mathcal{N}}(\bm{\mu}_{\text{train}}^i,t^j)\}_{i,j=1}^{{N_{\text{train}}},{N_t}}$.}
\Statex{{\bf Output:} Trained GCA-ROMs $\{\tilde{\mathbf{u}}_{j}^{\text{GCA-ROM}}(\boldsymbol{\mu})\}_{j=1}^{N_{\text{gca}}}$ and OpInf model $\tilde{\bm{g}}_{3}(t)$.}
\vskip+10pt
\State  Collect snapshots into the high-fidelity snapshot tensor $\mathbf{S}_{\text{hf}} \in \mathbb{R}^{N_{\text{train}}\times N_h\times N_t}$.
\State Perform a TT decomposition to approximate $\mathbf{S}_{\text{hf}}\approx \bm{G}_1 \otimes \bm{G}_2\otimes \bm{G}_3 $.
\State Train the $N_{\text{gca}}$ GCA-ROMs $\tilde{\mathbf{u}}_{j}^{\text{GCA-ROM}}(\boldsymbol{\mu})$ from the high-fidelity solutions $\{\mathbf{u}_{\mathcal{N}}(\bm{\mu}_{\text{train}}^i, t_{\text{gca}}^j)\}_{i,j=1}^{N_{\text{train}},N_{\text{gca}}}.$
\State Construct the operator inference model $\tilde{\bm{g}}_{3}(t)$ using data $\{\bm{g}_3^1, \bm{g}_3^2, \ldots, \bm{g}_3^{N_t}\}$ to predict the time-dependent core.
\end{algorithmic}

\vskip+10pt
\begin{algorithmic}[1]
\Statex{{\bf \underline{Online Phase}}}
\Statex{{\bf Input:} Parameter $\bm{\mu}^* \in \mathcal{P}$ and time $t^* \in [0, t_f]$ corresponding to the solution to be predicted.}
\Statex{{\bf Output:} Predicted solution $\tilde{\mathbf{u}}^{\text{LF}}(\bm{\mu}^{*}, t^{*})$.}
\vskip+10pt
\State Compute the GCA-ROM solutions $\tilde{\mathbf{u}}_j^{\text{GCA-ROM}}(\bm{\mu}^*)$ for $\bm{\mu}^* \in \mathcal{P}$ for $j = 1,\ldots, N_{\text{gca}}$, and build the predicted matrix  
\[
\mathbf{S}_{\text{gca}}(\bm{\mu}^*) = \begin{bmatrix} 
\tilde{\mathbf{u}}_{1}^{\text{GCA-ROM}}(\bm{\mu}^*)\ |\ \tilde{\mathbf{u}}_{2}^{\text{GCA-ROM}}(\bm{\mu}^*)\ |\ \cdots\ |\ \tilde{\mathbf{u}}_{N_{\text{gca}}}^{\text{GCA-ROM}}(\bm{\mu}^*)
\end{bmatrix}.
\]  
\State Solve the minimization problem:  
\begin{equation*}
	\tilde{\bm{g}}_1^{\text{LF}}(\bm{\mu}^*) = \underset{\hat{\bm{c}} \in \mathbb{R}^{r_1}}{\operatorname{argmin}} \left\| \mathbf{S}_{\text{gca}}(\bm{\mu}^*) - \hat{\bm{c}} \otimes \bm{G}_2 \otimes \bm{G}_3^{\text{gca}} \right\|^2.
\end{equation*}
\State Use the operator inference model to predict the core $\tilde{\bm{g}}_3(t^*)$.  
\State Compute $\tilde{\mathbf{u}}^{\text{LF}}(\bm{\mu}^*, t^*)$ by assembling low-fidelity coefficient $\tilde{\bm{g}}_{1}^{\text{LF}}(\bm{\mu}^*)$ and predicted core $\tilde{\bm{g}}_3(t^*)$ as:
\begin{equation*}
	\tilde{\mathbf{u}}^{\text{LF}}(\bm{\mu}^{*}, t^{*}) = \tilde{\bm{g}}_{1}^{\text{LF}}(\bm{\mu}^*)\otimes\bm{G}_{2}\otimes \tilde{\bm{g}}_{3}(t^*).
\end{equation*}
\end{algorithmic}
\end{algorithm}

The offline phase starts by collecting the parametric trajectories in time into the high-fidelity snapshot tensor $\mathbf{S}_{\text{hf}} \in \mathbb{R}^{N_{\text{train}}\times N_h\times N_t}$, and performing a TT decomposition to obtain its low-rank approximation as $\mathbf{S}_{\text{hf}} \approx \bm{G}_1 \otimes \bm{G}_2\otimes \bm{G}_3$. 
Then, we aim at decoupling the temporal and parametric dependence, by selecting a suitably small number $N_{\text{gca}}\leq N_t$ of time points $\{t_{\text{gca}}^j\}_{j=1}^{N_{\text{gca}}}$, and training a GCA-ROM for each of them from the dataset $\{\bm{\mu}_{\text{train}}^i, \mathbf{u}_{\mathcal{N}}(\bm{\mu}_{\text{train}}^i, t_{\text{gca}}^j)\}_{i=1}^{N_{\text{train}}}$, where $j = 1, 2, \ldots, N_{\text{gca}}$. 
This way, we exploit a geometrically consistent nonlinear reduction for the behavior of the system at fixed time instances, discarding the time complexity but capturing the parametric influence.
Finally, an OpInf surrogate model $\tilde{\bm{g}}_{3}(t)$ is built from the time coefficients $\{\bm{g}_3^1, \bm{g}_3^2, \ldots, \bm{g}_3^{N_t}\}$ to predict the evolution of the time-dependent core.

At inference time, i.e.\ during the online phase, we aim at predicting and extrapolating a low-fidelity approximation of the trajectory corresponding to a new parameter $\bm{\mu}^* \in \mathcal{P}$ at any given time $t^*\in[0,t_f]$ as 
\begin{equation}\label{eq:lf_approx}
\tilde{\mathbf{u}}^{\text{LF}}(\bm{\mu}^{*}, t^{*}) = \tilde{\bm{g}}_{1}^{\text{LF}}(\bm{\mu}^*)\otimes\bm{G}_{2}\otimes \tilde{\bm{g}}_{3}(t^*).
\end{equation}
Toward this goal, we start by evaluating the $N_{\text{gca}}$ trained GCA-ROM models, and we obtain the low-dimensional representation of the solution snapshots denoted by 
\[
\mathbf{S}_{\text{gca}}(\bm{\mu}^*) = \begin{bmatrix} 
\tilde{\mathbf{u}}_{1}^{\text{GCA-ROM}}(\bm{\mu}^*)\ |\ \tilde{\mathbf{u}}_{2}^{\text{GCA-ROM}}(\bm{\mu}^*)\ |\ \cdots\ |\ \tilde{\mathbf{u}}_{N_{\text{gca}}}^{\text{GCA-ROM}}(\bm{\mu}^*)
\end{bmatrix}.
\]  
Then, an approximation $\tilde{\bm{g}}_1^{\text{LF}}(\bm{\mu}^*)$ of the parametric core $ \bm{g}_{1}(\bm{\mu}^*)$ can be obtained by seeking the coefficients $\hat{\bm{c}} \in \mathbb{R}^{r_1}$ that fit $\mathbf{S}_{\text{gca}}(\bm{\mu}^*)$ in the TT format at $N_{\text{gca}}$ time instances, i.e.\ by solving the following minimization problem:
\begin{equation}
	\label{DMDLS}
	\tilde{\bm{g}}_1^{\text{LF}}(\bm{\mu}^*) = \underset{\hat{\bm{c}} \in \mathbb{R}^{r_1}}{\operatorname{argmin}} \left\| \mathbf{S}_{\text{gca}}(\bm{\mu}^*) - \hat{\bm{c}} \otimes \bm{G}_2 \otimes \bm{G}_3^{\text{gca}} \right\|^2,
\end{equation}
where $\bm{G}_3^{\text{gca}} = \left[\bm{g}_{3}(t_{\text{gca}}^1), \bm{g}_{3}(t_{\text{gca}}^2), \cdots, \bm{g}_{3}(t_{\text{gca}}^{N_{\text{gca}}})\right]$.

By evaluating the OpInf model to predict the core $\tilde{\bm{g}}_3(t^*)$, we finally have all the ingredients to reconstruct the solution at time $t^*\in[0,t_f]$ for any given parameter $\bm{\mu}^*\in \mathcal{P}$ as $\tilde{\mathbf{u}}^{\text{LF}}(\bm{\mu}^{*}, t^{*}) = \tilde{\bm{g}}_{1}^{\text{LF}}(\bm{\mu}^*)\otimes\bm{G}_{2}\otimes \tilde{\bm{g}}_{3}(t^*)$.

\subsection{Multi-fidelity extension of LF-TTOI via DeepONet}
\begin{figure}[htbp]
\begin{center}
    \includegraphics[width=0.9\textwidth]{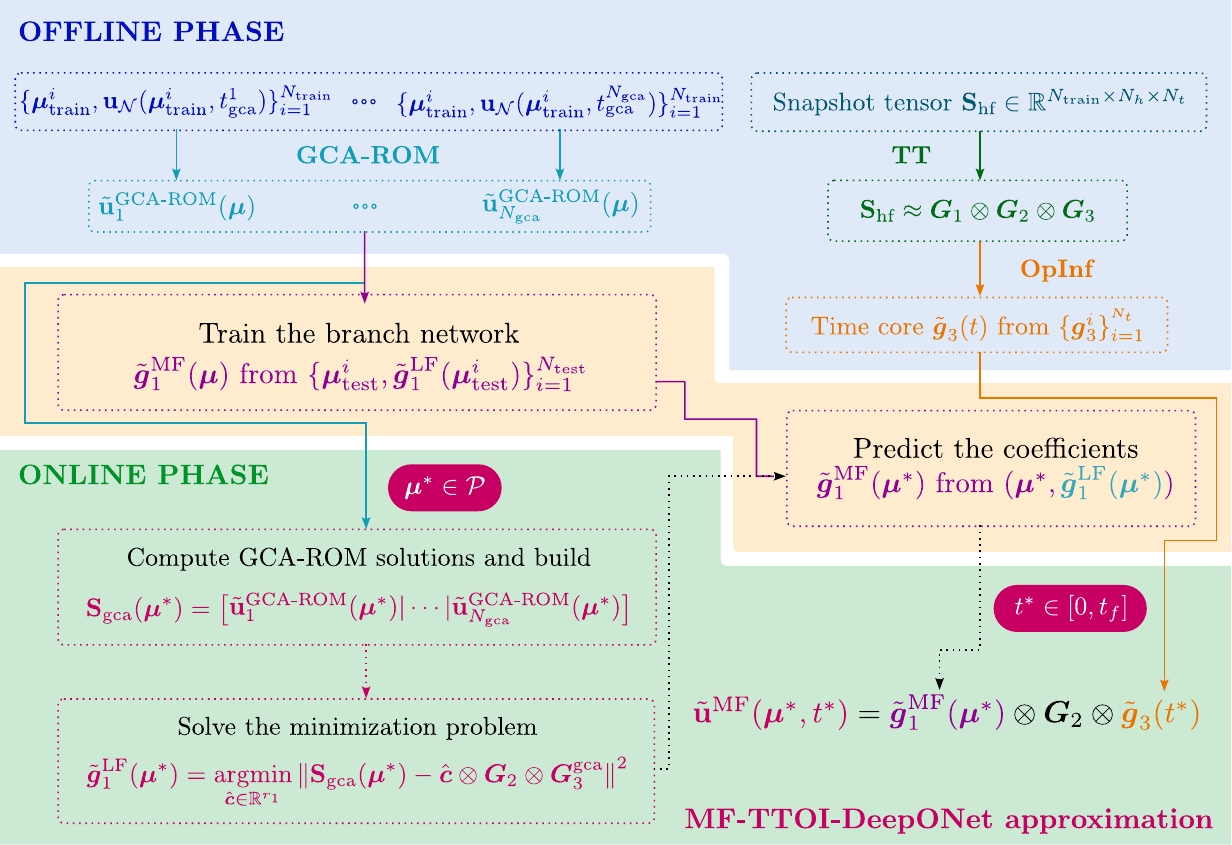}
\end{center}
\caption{The framework of MF-TTOI-DeepONet.}
\label{mf}
\end{figure}

LF-TTOI relies on the GCA-ROM surrogates for different time slices to compute the parametric coefficients $\tilde{\bm{g}}_1^{\text{LF}}(\bm{\mu})$. When the high-fidelity dynamics is particularly complex and only scarce data are available, the accuracy of the approach could be affected. This could be a potential limitation of the described approach, e.g.\ for benchmarks where the physical/geometrical parametrization changes the inherent velocity of the dynamics, and more information are needed to recover the global behavior. We remark that also the minimization problem in Equation \eqref{DMDLS} introduces an additional source of error in the process that need to be taken into account. 
For these reasons, we propose an advanced strategy from the perspective of Neural Operator, namely DeepONet \cite{Lu2019LearningNO}.
In particular, we developed a multi-fidelity strategy, denoted as MF-TTOI-DeepONet, where \( \tilde{\bm{g}}_1^{\text{LF}}(\bm{\mu}) \) is regarded as the low-fidelity approximation of the parameterized coefficients ${\bm{g}}_1(\bm{\mu})$, that we use to improve the performance of the model by capturing its discrepancies from the original coefficients. This allows us to generalize across all parameters, and ultimately achieving more accurate predictions.
The proposed multi-fidelity strategy MF-TTOI-DeepONet is detailed in Algorithm \ref{MF-TTOI-DeepONet}, while a schematic representation of the framework is shown in Figure \ref{mf}. 

\begin{algorithm}[t]
\caption{MF-TTOI-DeepONet}
\label{MF-TTOI-DeepONet}
\vskip+10pt
\begin{algorithmic}[1]
\Statex{{\bf \underline{Offline Phase}} }
\Statex{{\bf Input:} Training dataset with parameter samples $\Xi_{\mathrm{train}}=\{\bm{\mu}_{\text{train}}^i\}_{i=1}^{N_{\text{train}}}$, time steps $\{t^j\}_{j=1}^{N_t}$,}
\LineComment{GCA time points $\{t_{\text{gca}}^j\}_{j=1}^{N_{\text{gca}}}$, and high-fidelity solutions $\{\mathbf{u}_{\mathcal{N}}(\bm{\mu}_{\text{train}}^i,t^j)\}_{i,j=1}^{{N_{\text{train}}},{N_t}}$.}
\Statex{{\bf Output:} Trained GCA-ROMs $\{\tilde{\mathbf{u}}_{j}^{\text{GCA-ROM}}(\boldsymbol{\mu})\}_{j=1}^{N_{\text{gca}}}$, operator inference model $\tilde{\bm{g}}_{3}(t)$,}
\LineComment{ and trained branch net $\tilde{\bm{g}}_1^{\text{MF}}(\bm{\mu})$.}
\vskip+10pt
\State  Collect snapshots into the high-fidelity snapshot tensor $\mathbf{S}_{\text{hf}} \in \mathbb{R}^{N_{\text{train}}\times N_h\times N_t}$.
\State Perform a TT decomposition to approximate $\mathbf{S}_{\text{hf}}\approx \bm{G}_1 \otimes \bm{G}_2\otimes \bm{G}_3 $.
\State Train the $N_{\text{gca}}$ GCA-ROMs $\tilde{\mathbf{u}}_{j}^{\text{GCA-ROM}}(\boldsymbol{\mu})$ from the high-fidelity solutions $\{\mathbf{u}_{\mathcal{N}}(\bm{\mu}_{\text{train}}^i, t_{\text{gca}}^j)\}_{i,j=1}^{N_{\text{train}},N_{\text{gca}}}.$
\State Construct the operator inference model $\tilde{\bm{g}}_{3}(t)$ using data $\{\bm{g}_3^1, \bm{g}_3^2, \ldots, \bm{g}_3^{N_t}\}$ to predict the time-dependent core.

\State{Train the branch network $\tilde{\bm{g}}_{1}^{\text{MF}}$ with input $\bm{\mu}$ and low-fidelity coefficients $\tilde{\bm{g}}_1^{\text{LF}}(\bm{\mu})$, based on low-fidelity GCA-ROM data $\{\tilde{\mathbf{u}}^{\text{GCA-ROM}}(\bm{\mu}_{\text{test}}^i)\}_{i,j=1}^{N_{\text{test}},N_{\text{gca}}}$, using the loss function Equation (\ref{mfloss1}) to enhance generalization.}

\State{Fine-tune the branch network $\tilde{\bm{g}}_1^{\text{MF}}$ with input $\bm{\mu}$ and low-fidelity coefficients $\tilde{\bm{g}}_1^{\text{LF}}(\bm{\mu})$, using the high-fidelity solutions $\{\mathbf{u}_{\mathcal{N}}(\bm{\mu}_{\text{train}}^i,t^j)\}_{i,j=1}^{N_{\text{train}},N_t}$, with the loss function Equation (\ref{mfloss}).}

\end{algorithmic}

\vskip+10pt
\begin{algorithmic}[1]
\Statex{{\bf \underline{Online Phase}}}
\Statex{{\bf Input:} Parameter $\bm{\mu}^* \in \mathcal{P}$ and time $t^* \in [0, t_f]$ corresponding to the solution to be predicted.}
\Statex{{\bf Output:} Predicted solution $\tilde{\mathbf{u}}^{\text{MF}}(\bm{\mu}^{*}, t^{*})$.}
\vskip+10pt
\State Compute the GCA-ROM solutions $\tilde{\mathbf{u}}_j^{\text{GCA-ROM}}(\bm{\mu}^*)$ for $\bm{\mu}^* \in \mathcal{P}$ for $j = 1,\ldots, N_{\text{gca}}$, and build the predicted matrix  
\[
\mathbf{S}_{\text{gca}}(\bm{\mu}^*) = \begin{bmatrix} 
\tilde{\mathbf{u}}_{1}^{\text{GCA-ROM}}(\bm{\mu}^*)\ |\ \tilde{\mathbf{u}}_{2}^{\text{GCA-ROM}}(\bm{\mu}^*)\ |\ \cdots\ |\ \tilde{\mathbf{u}}_{N_{\text{gca}}}^{\text{GCA-ROM}}(\bm{\mu}^*)
\end{bmatrix}.
\]  
\State Solve the minimization problem:  
\begin{equation*}
	\tilde{\bm{g}}_1^{\text{LF}}(\bm{\mu}^*) = \underset{\hat{\bm{c}} \in \mathbb{R}^{r_1}}{\operatorname{argmin}} \left\| \mathbf{S}_{\text{gca}}(\bm{\mu}^*) - \hat{\bm{c}} \otimes \bm{G}_2 \otimes \bm{G}_3^{\text{gca}} \right\|^2.
\end{equation*}
\State Use the branch net 
to predict the  multi-fidelity coefficient $\tilde{\bm{g}}_1^{\text{MF}}(\bm{\mu}^*)$ with input $(\bm{\mu}^*, \tilde{\bm{g}}_1^{\text{LF}}(\bm{\mu}^*))$.
\State Use the operator inference model to predict the core $\tilde{\bm{g}}_3(t^*)$. 
\State Compute $\tilde{\mathbf{u}}^{\text{MF}}(\bm{\mu}^*, t^*)$ by assembling $\tilde{\bm{g}}_1^{\text{MF}}(\bm{\mu}^*)$, and predicted core $\tilde{\bm{g}}_3(t^*)$ as:
\begin{equation*}
	\tilde{\mathbf{u}}^{\text{MF}}(\bm{\mu}^{*}, t^{*}) = \tilde{\bm{g}}_{1}^{\text{MF}}(\bm{\mu}^*)\otimes\bm{G}_{2}\otimes \tilde{\bm{g}}_{3}(t^*).
\end{equation*}
\end{algorithmic}
\end{algorithm}

In MF-TTOI-DeepONet, the multi-fidelity approximation of the original snapshot corresponding to a new parameter $\bm{\mu}^* \in \mathcal{P}$ at any given time $t^*\in[0,t_f]$ is sought as:
\begin{equation}\label{eq:mf_approx}
	\tilde{\mathbf{u}}^{\text{MF}}(\bm{\mu}^{*}, t^{*}) = \tilde{\bm{g}}_{1}^{\text{MF}}(\bm{\mu}^*)\otimes\bm{G}_{2}\otimes \tilde{\bm{g}}_{3}(t^{*}),
\end{equation}
where, in analogy with the terminology of DeepONet, the branch network is used to obtain a multi-fidelity approximation $\tilde{\bm{g}}_{1}^{\text{MF}}(\bm{\mu})$ of the parameter-dependent core, while the term $\bm{G}_{2} \otimes \tilde{\bm{g}}_{3}(t)$ is considered to be the outputs of the trunk network. 
Specifically, the branch network providing a multi-fidelity estimation of ${\bm{g}}_1(\bm{\mu})$, is represented as the following correction of the LF-TTOI approximation:
$$
 \tilde{\bm{g}}_{1}^{\text{MF}}(\bm{\mu}, \tilde{\bm{g}}_1^{\text{LF}}(\bm{\mu})) = \tilde{\bm{g}}_1^{\text{LF}}(\bm{\mu}) + f_{nn}(\bm{\mu}),
$$
where $f_{nn}(\bm{\mu})$ is representing the correction, and it is modelled as a FNN with the parameter $\bm{\mu}$ as input.

To further leverage the predictions of the GCA-ROMs and fully utilize the collected high-fidelity data, we divide the training of the branch network $\tilde{\bm{g}}_{1}^{\text{MF}}$ into two stages.  In the first stage, a pre-training strategy is employed to enhance the diversity of training samples and improve the network's generalization capability. Here, the network is trained on the low-fidelity dataset $\{\tilde{\mathbf{u}}_j^{\text{GCA-ROM}}(\bm{\mu}_{\text{test}}^i)\}_{i,j=1}^{N_{\text{test}},N_{\text{gca}}}$, generated by the well-trained GCA-ROM, along with the corresponding low-fidelity coefficients for the parameter sets $\Xi_{\mathrm{test}}=\{\bm{\mu}_{\text{test}}^i\}_{i=1}^{N_{\text{test}}}$. The loss function is defined as:  
\begin{equation}
\label{mfloss1}
\mathcal{L}_{\mathrm{LF}} = \frac{1}{N_{\text{test}}} \sum_{i=1}^{N_{\text{test}}} \frac{\left\| \tilde{\mathbf{u}}^{\text{MF}}\left(\bm{\mu}_{\text{test}}^i, \bm{T}_{\text{gca}}\right) - \tilde{\mathbf{u}}^{\text{GCA-ROM}}\left(\bm{\mu}_{\text{test}}^i, \bm{T}_{\text{gca}}\right) \right\|_2^2}{\left\| \tilde{\mathbf{u}}^{\text{GCA-ROM}}\left(\bm{\mu}_{\text{test}}^i, \bm{T}_{\text{gca}}\right) \right\|_2^2},
\end{equation} 
where $\bm{T}_{\text{gca}} = [t_{\text{gca}}^1, \ldots, t_{\text{gca}}^{N_{\text{gca}}}]$ denotes the time instances associated with the GCA-ROM predictions.  

In the second stage, the branch network $\tilde{\bm{g}}_{1}^{\text{MF}}$ is further fine-tuned using the high-fidelity dataset to improve prediction accuracy. The corresponding loss function is given by:  
\begin{equation}
\label{mfloss}
\mathcal{L}_{\mathrm{MF}} = \frac{1}{N_{\text{train}}} \sum_{i=1}^{N_{\text{train}}} \frac{\left\| \tilde{\mathbf{u}}^{\text{MF}}\left(\bm{\mu}_{\text{train}}^i, \bm{T}_{\text{train}}\right) - \mathbf{u}_{\mathcal{N}}\left(\bm{\mu}_{\text{train}}^i, \bm{T}_{\text{train}}\right) \right\|_2^2}{\left\| \mathbf{u}_{\mathcal{N}}\left(\bm{\mu}_{\text{train}}^i, \bm{T}_{\text{train}}\right) \right\|_2^2},
\end{equation}  
where $\bm{T}_{\text{train}} = [t^1, \ldots, t^{N_t}]$ represents the time steps in the training dataset.  
Thus, the collected high-fidelity data is not only used for the TT decomposition but also reused to optimize the branch network during training.  

Compared to LF-TTOI, achieving more reliable accuracy requires training an additional operator network during the offline phase. However, it is important to note that no extra high-fidelity data is needed for this process.

\subsection{Computational costs}
Here we discuss the computational cost of the ingredients of the proposed approaches. Let us start with the analysis of Algorithm \ref{LF-TTOI}, ignoring the cost for the training and evaluation of the neural network, and focusing on the scalability w.r.t.\ the characteristic dimensions of the different quantities. In the offline stage, the cost of performing TT decomposition on the snapshot tensor $\mathbf{S}_{\text{hf}}$ is $\mathcal{O}(N_{\text{train}}^2 N_h N_t + \hat{r}^2N_h^{2}N_t)$,  where $\hat{r} = \max\{r_1, r_2\}$, while the cost of building the OpInf model is $\mathcal{O}(N_{t} {r}_2\bar{n}^3)$. In the online stage, for any given parameter $\bm{\mu}^* \in \mathcal{P}$, the cost of obtaining the low-fidelity coefficients $\tilde{\bm{g}}_1^{\text{LF}}(\bm{\mu}^*)$ via the optimization problem in Equation \eqref{DMDLS} is  given by $\mathcal{O}(N_{h}{N_{\text{gca}}}r_1^3)$. The cost of solving the least-squares problem in Equation \eqref{ls} grows exponentially with the complexity of the nonlinearity. Therefore, in this work, the OpInf model in Equation \eqref{op} is limited to polynomial nonlinear terms of order no greater than $L=2$.  As an example, by taking a second-order nonlinear model, the OpInf cost is $\mathcal{O}({r}_2n_2)$ for each step, while the cost for the high-fidelity reconstruction $\tilde{\mathbf{u}}^{\text{LF}}(\bm{\mu}^{*}, t^{*}) $ is $\mathcal{O}(N_{h}\hat{r}^2)$.

As concerns the multi-fidelity strategy presented in Algorithm \ref{MF-TTOI-DeepONet}, during the offline stage, it requires an additional step to collect low-fidelity coefficients $\{\tilde{\bm{g}}_1^{\text{LF}}(\bm{\mu}_{\text{train}}^i)\}_{i=1}^{N_{\text{train}}}$, whose computational cost is given by $\mathcal{O}(N_{\text{train}}N_{h}{N_{\text{gca}}}r_1^3)$, and train the branch network $\tilde{\bm{g}}_1^{\text{MF}}(\bm{\mu})$.

\section{Numerical results}
\label{Numerical examples}
In this section, we analyze the performance of the proposed methods on three benchmark problems: heat conduction in a thermal block, advection-diffusion problem in a square domain with a square obstacle, and the Navier-Stokes system, and compare the proposed method with the performance of the original GCA-ROM approach and several state-of-the-art techniques. The datasets are generated using the open-source tools FEniCS \cite{alnaes2015fenics} and RBniCS \cite{RozzaRealTimeReduced2024}.
For the error analysis, we define the relative error of the different methods for a given parameter \(\boldsymbol{\mu} \in \mathcal{P}\) and time $t\in[0,t_f]$ as  
\[
\varepsilon_{\square}(\boldsymbol{\mu},t) = \frac{\|\mathbf{u}_{\mathcal{N}}(\boldsymbol{\mu},t) - \tilde{\mathbf{u}}^{\square}(\boldsymbol{\mu},t)\|_2}{\|\mathbf{u}_{\mathcal{N}}(\boldsymbol{\mu},t)\|_2},
\]
the mean relative error over the testing parameter set \( \Xi_{\mathrm{test}} \) as  
\[
\bar{\varepsilon}_{\square}(t) = \frac{1}{N_{\text{test}}}  \sum_{i=1}^{N_{\text{test}}}  \varepsilon_{\square}(\boldsymbol{\mu}_{\text{test}}^i,t),
\]
the mean relative error over the time domain \(\boldsymbol{T}_{\text{test}} = \{t_{\text{test}}^1,\cdots,t_{\text{test}}^{N_{s}}\} \in [0,t_f] \) as  
\[
\bar{\varepsilon}_{\square}(\bm{\mu}) = \frac{1}{N_{s}} \sum_{i=1}^{N_{s}}  \varepsilon_{\square}(\boldsymbol{\mu},t_{\text{test}}^i),
\]
and the global mean relative error over the test parameter-time domain as  
\[
\bar{\bar{\varepsilon}}_{\square} = \frac{1}{ N_s} \sum_{i=1}^{N_s} \bar{ \varepsilon}_{\square}(t_{\text{test}}^i),
\]
where \( \square \in \{ \text{MF}, \text{LF} \} \),  
and \( \tilde{\mathbf{u}}^{\square}(\boldsymbol{\mu}, t) \) represents the solution predicted by the corresponding model MF-TTOI-DeepONet or LF-TTOI.
The TT decomposition for all the benchmarks is performed with a tolerance \( \varepsilon_\text{TT} = 10^{-8} \).  Table~\ref{tvalue} shows the time intervals $\boldsymbol{T}_{\text{train}}$ and $\boldsymbol{T}_{\text{test}}$, which are uniformly discretized with the specified step sizes, and the selected time instances $\boldsymbol{T}_{\text{gca}}$ for the GCA-ROM method for each benchmark problem.
\begin{table}[h!]
\centering
\begin{tabular}{@{}ccccc@{}}
\toprule
\textbf{Application}      &  step size &  $\boldsymbol{T}_{\text{train}}$  &  $\boldsymbol{T}_{\text{test}}$    & $\boldsymbol{T}_{\text{gca}}$        \\ \midrule
Heat conduction  & 0.05& [0, 1]&[0, 3] &[0.05, 0.2, 0.4, 0.6, 0.8, 1]\\ 
Advection-diffusion & 0.02  & [0, 1.5]&[0, 2] &[0.02, 0.38, 0.76, 1.12, 1.5] \\ 
Navier-Stokes cylinder &0.01& [0, 4.49]&[0, 6] &[0.01, 2.25,  4.49] \\ 
Navier-Stokes rectangle & 0.01 & [0, 0.75]&[0, 1] &[0.01, 0.13, 0.25, 0.38, 0.5, 0.62, 0.75]\\\bottomrule
\end{tabular}
\caption{The information of $\boldsymbol{T}_{\text{train}}$, $\boldsymbol{T}_{\text{gca}}$ and $\boldsymbol{T}_{\text{test}}$ for the different benchmarks.}\label{tvalue}
\end{table}

\subsection{Heat conduction in a thermal block}
\label{heateq}
To illustrate the proposed methods,  we first consider a time-dependent heat conduction problem in a square domain $\Omega=[-1,1]^2$ in Figure \ref{heatdomain}. 

\definecolor{bluemathlab}{HTML}{065895}
\definecolor{greenmathlab}{HTML}{13ad8e}
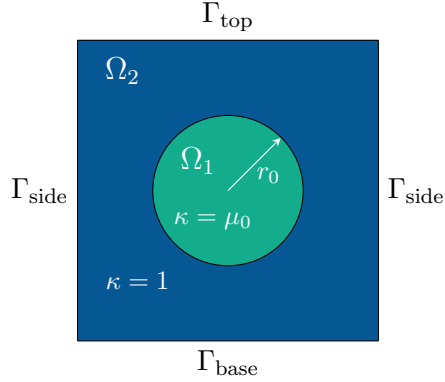
\begin{figure}[htb]
\begin{center}
            \centering
        \begin{tikzpicture}[scale=2]
        \fill [bluemathlab] (-1,-1) rectangle (1,1);
        \fill [greenmathlab] (0,0) circle (0.5cm);
        \draw[-] (0,0) circle (0.5);
        \draw (-1,-1) -- (1,-1) node[pos=0.5,sloped,below] {\large{$\Gamma_{\textrm{base}}$}};
        \draw (1,-1) -- (1,1) node[pos=0.5,right] {\large{$\Gamma_{\textrm{side}}$}};
        \draw (1,1) -- (-1,1) node[pos=0.5,sloped,above] {\large{$\Gamma_{\textrm{top}}$}};
        \draw (-1,1) -- (-1,-1) node[pos=0.5,left] {\large{$\Gamma_{\textrm{side}}$}};
        \node[white] at (-0.20,0.2) {\Large{$\Omega_1$}};
        \node[white] at (-0.7,0.8) {\Large{$\Omega_2$}};
        \node[white] at (-0.6,-0.6) {\normalsize{$\kappa = 1$}};
        \node[white] at (-0.1,-0.2) {\normalsize{$\kappa = \mu_0$}};
        \draw[white,->,>=stealth] (0,0) -- (0.35,0.35);
        \node[white] at (0.27,0.1) {\normalsize{$r_0$}};
        \end{tikzpicture}
\end{center}
\caption{Geometrical setup of the heat conductivity problem.}
\label{heatdomain}
\end{figure}

The domain is divided into two subdomains, a disk $\Omega_1$ centered at the origin with radius $r_0=0.5$, and its surrounding area $\Omega_2 = \Omega/{\Omega_1}$. The temperature field \( u(t; \bm{\mu}) \) is governed by a parameterized parabolic PDE, where the parameter vector \( \bm{\mu} = (\mu_0, \mu_1) \in \mathbb{P} = [0.1,10]\times[-1,1]\) contains the thermal conductivity $\kappa$ in $\Omega$ and the heat flux over the bottom boundary $\Gamma_\text{base} $.

The governing equation, along with boundary and initial conditions, is expressed as:
$$
\begin{cases}
    \frac{\partial u}{\partial t} - \nabla \cdot \big(\kappa_{\mu_0} \nabla u \big) = 0, &\text{in } \Omega \times (0, t_f], \\
    u(x,0; \bm{\mu}) = 0, &\text{in } \Omega, \\
    u(x,t; \bm{\mu}) = 0, &\text{on } \Gamma_{\text{top}} \times (0, t_f], \\
    \kappa_{\mu_0} \nabla u \cdot \mathbf{n} = 0, &\text{on } \Gamma_{\text{side}} \times (0, t_f], \\
    \kappa_{\mu_0} \nabla u \cdot \mathbf{n} = \mu_1, &\text{on } \Gamma_{\text{base}} \times (0, t_f],
\end{cases}
$$
where \( \kappa_{\mu_0} \) represents the thermal conductivity, defined as \( \kappa_{\mu_0} = \mathbf{1}_{\Omega_2} + \mu_0 \mathbf{1}_{\Omega_1} \), where \( \mathbf{1}_{\Omega_1} \) and \( \mathbf{1}_{\Omega_2} \) are the indicator functions for the two subdomains, and $\mathbf{n}$ denotes the outer normal to the boundaries $ \Gamma_{\text{side}} $ and $\Gamma_{\text{base}}$.

To generate the high-fidelity snapshots, we utilize the  FE method with $\mathbb{P}_1$ elements, resulting in $N_h = 304$ degrees of freedom. We employ the backward Euler method with a time step size $\Delta t = 0.05$ to simulate the dynamics up to $t_f = 3$.  Finally, the parametric dataset \( \Xi \) consists of  \( N_{\bm{\mu}} = 100 \) parameters corresponding to 10 equispaced parameter samples in each dimension. By setting $r_t=50\%$, a set of $N_\text{train} = 50$ randomly selected parameter values \( \Xi_\text{train} \) and the corresponding solutions are used for training, while the remaining snapshots are used to build the test dataset  \( \Xi_\text{test} \).

\begin{figure}[htbp]
\begin{center}
    \includegraphics[width=0.45\textwidth]{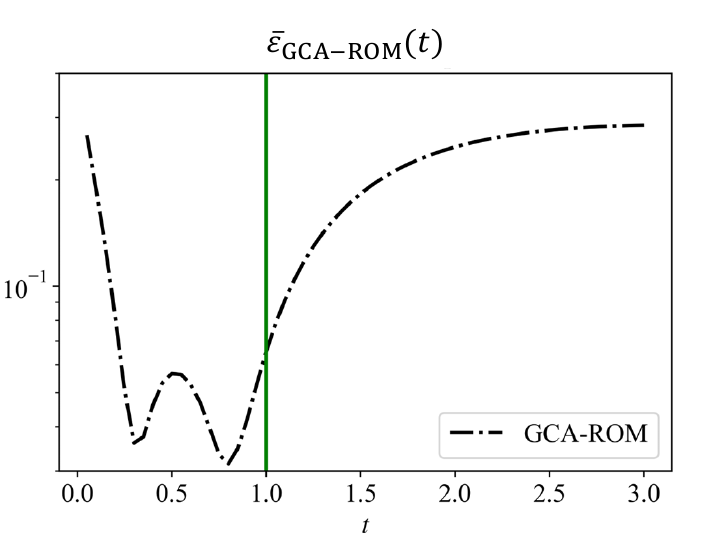}
\end{center}
\vskip-20pt
\caption{Error ${\bar{\varepsilon}}_{\text{GCA-ROM}}(t)$ of the plain GCA-ROM on test parameter set.}
\label{fig:gca}
\end{figure}
\begin{figure}[htbp]
\begin{center}
    \includegraphics[width=0.8\textwidth]{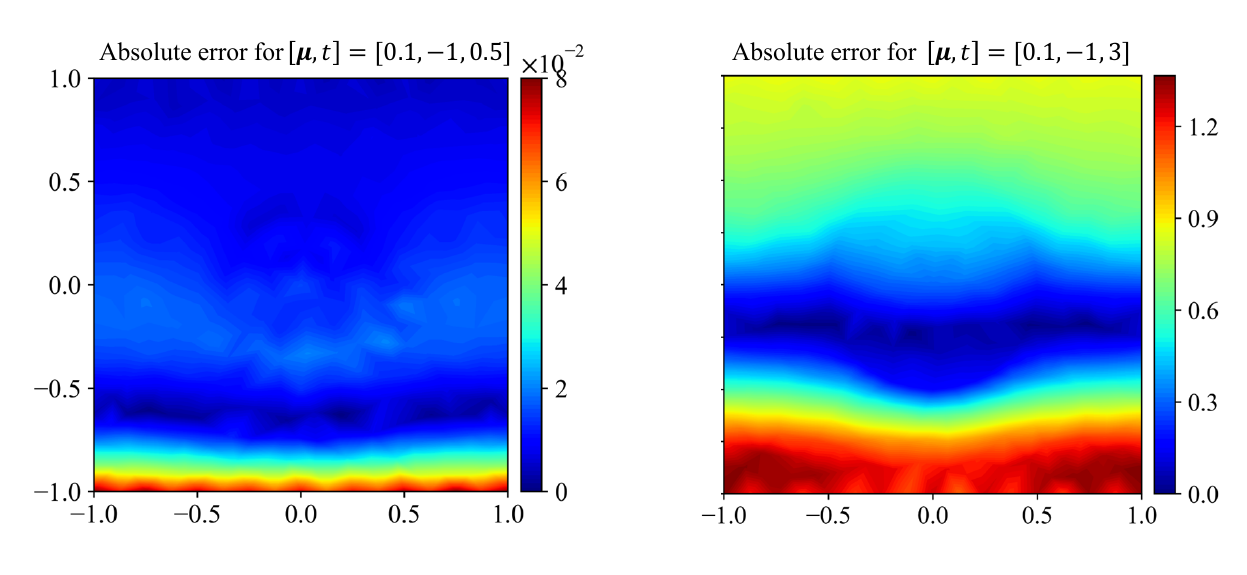}
\end{center}
\vskip-20pt
\caption{Absolute errors of the plain GCA-ROM at time instances $t = 0.5$ and $t = 3$ for $\bm{\mu} = [0.1, -1]$.}
\label{fig:gcaabs}
\end{figure}

We begin by illustrating the performance of GCA-ROM on time-dependent problems, where \( t \) is treated as an additional parameter. Figure \ref{fig:gca} shows the error ${\bar{\varepsilon}}_{\text{GCA-ROM}}(t)$ of GCA-ROM  on the test parameter set. It can be observed that, as expected, the errors increase rapidly outside the training time range \([0, 1]\), thus exhibiting bad prediction capabilities in the extrapolation regime due to the lack of temporal causality. In Figure \ref{fig:gcaabs}, we depict the absolute error of GCA-ROM at different time instances $[\bm{\mu}, t] = [0.1, -1, 0.5]$ and $[\bm{\mu}, t] = [0.1, -1, 3]$, from where it can be seen that the behavior of the error. This corroborates the thesis that, considering time as a simple parameter within these architectures, only allows to learn already observed patterns and correlations in the dataset, without really learning the underlying dynamics and the existing causality between subsequent snapshots. Moreover, we remark that both tricky initial conditions and fast changes in time of the dynamics could potentially jeopardize the learned surrogate model, as it can be seen already with this simple benchmark even within the training time range.

\begin{figure}[htbp]
\begin{center}
    \includegraphics[width=0.99\textwidth]{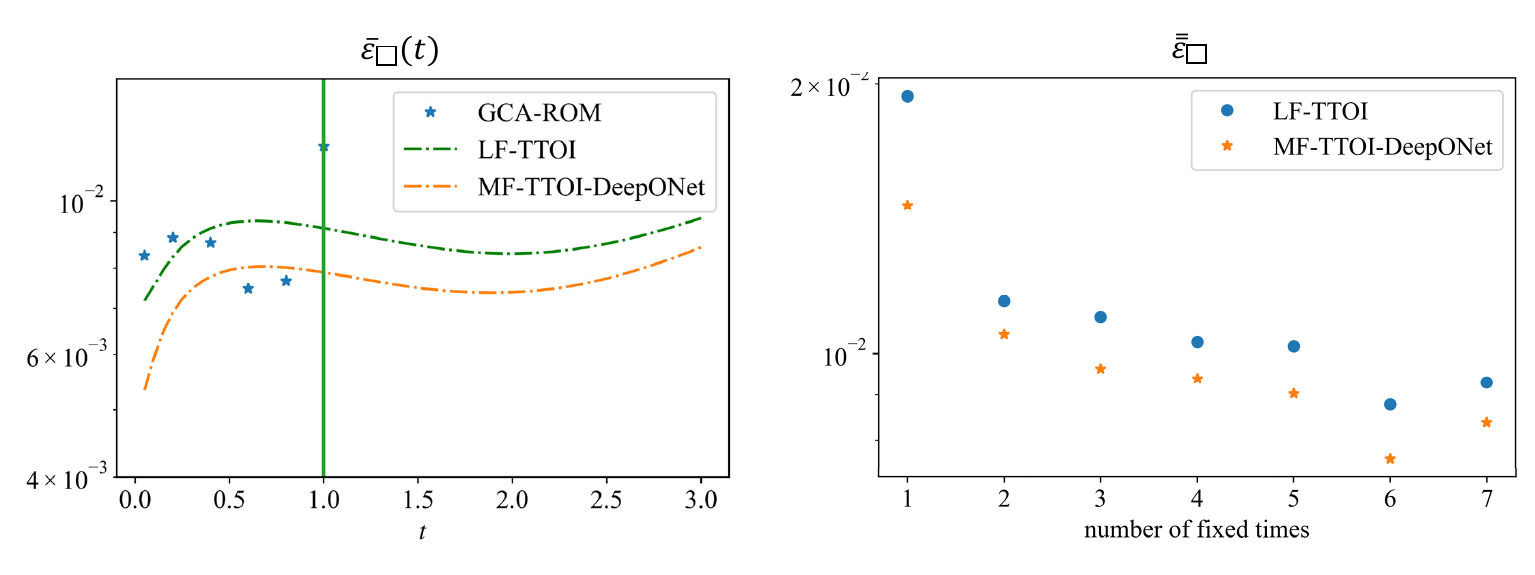}
\end{center}
\vskip-20pt
\caption{Heat conduction problem. Left:  Mean relative error ${\bar{\varepsilon}}_{\square}(t)$ of different methods on test parameters over the entire time domain. The vertical line at \( t = 1\) denotes the extrapolation regime. Right: Mean relative error $\bar{\bar{\varepsilon}}_{\square}$ when exploiting a different number of GCA-ROMs $N_{\text{gca}}$.}
\label{fig:mreheat1}
\end{figure}

To overcome the limitations of the plain GCA-ROM, the introduced methods exploit OpInf to learn the dynamics. Using the same training data, truncated via TT decomposition  with ranks $r_1=7$ and $r_2=11$ for parameter and time respectively, we evaluate the effectiveness of both LF-TTOI and MF-TTOI-DeepONet approaches. The left plot in Figure \ref{fig:mreheat1} shows the mean relative error ${\bar{\varepsilon}}_{\square}(t)$ on test parameter space \( \Xi_\text{test} \) and a comparison with the different techniques when choosing $N_{\text{gca}}=6$. 
For both proposed approaches, we observed a completely different behavior of the errors than the plain GCA-ROM approach.
It can be seen that, by decoupling the temporal and parametric influence, the LF-TTOI and MF-TTOI-DeepONet methodologies are able to obtain good approximations for both the training and extrapolation time ranges, providing accurate predictions with relative errors on the order of $10^{-3}$. 

We now test the effect of the number of GCA-ROM surrogates in time on the accuracy of proposed methods.
The right plot of Figure \ref{fig:mreheat1} shows the variation of the errors $\bar{\bar{\varepsilon}}_{\text{LF}}$  and $\bar{\bar{\varepsilon}}_{\text{MF}}$ w.r.t.\ the number of GCA-ROMs $N_{\text{gca}}$.
We start by noticing that increasing the number of GCA-ROMs helps in reducing the error as expected, until reaching a certain threshold, that is $N_{\text{gca}}>6$, after which we observe a saturation of the error itself indicating that $\bar{\bar{\varepsilon}}_{\text{LF}}$ and $\bar{\bar{\varepsilon}}_{\text{MF}}$ are not benefiting from the additional information, i.e.\ the dynamics is correctly resolved while more knowledge about the parametric variability could be needed. 

Surprisingly, even with two GCA-ROMs, the LF-TTOI is capable of achieving satisfactory precision, exploiting the spatio-temporal core to enhance the performance in comparison to GCA-ROM. 
The MF-TTOI-DeepONet, which demonstrates optimal performance across all time points, indicates that the integration of a branch network further enhances predictive accuracy.
This is attributed to the fact that, compared to LF-TTOI, MF-TTOI-DeepONet incorporates a FNN to model the residual between the low-fidelity parameter core computed by LF-TTOI and the high-fidelity parameter core. This residual learning strategy leverages high-fidelity snapshots to compensate for the limitations of LF-TTOI and improve the overall accuracy.

\begin{figure}[htb]
\begin{center}
    \includegraphics[width=0.99\textwidth]{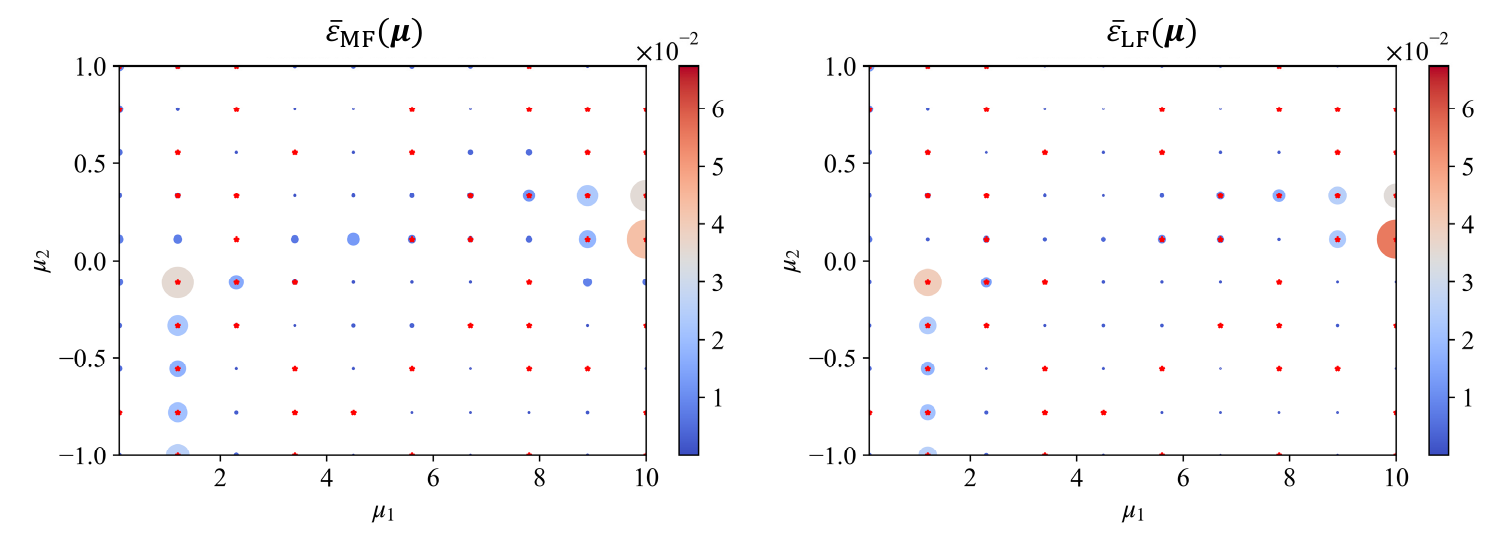}
\end{center}
\vskip-20pt
\caption{Heat conduction problem. Mean relative errors ${\bar{\varepsilon}}_{\mathrm{MF}}(\bm{\mu})$ for MF-TTOI-DeepONet (left) and ${\bar{\varepsilon}}_{\mathrm{LF}}(\bm{\mu})$ for LF-TTOI (right) on the dataset $\Xi$. Red markers corresponds to the parameters in the testing set $\Xi_{\text{test}}$.}
\label{fig:mreheat2}
\end{figure}

\begin{figure}[htb]
\begin{center}
    \includegraphics[width=0.99\textwidth]{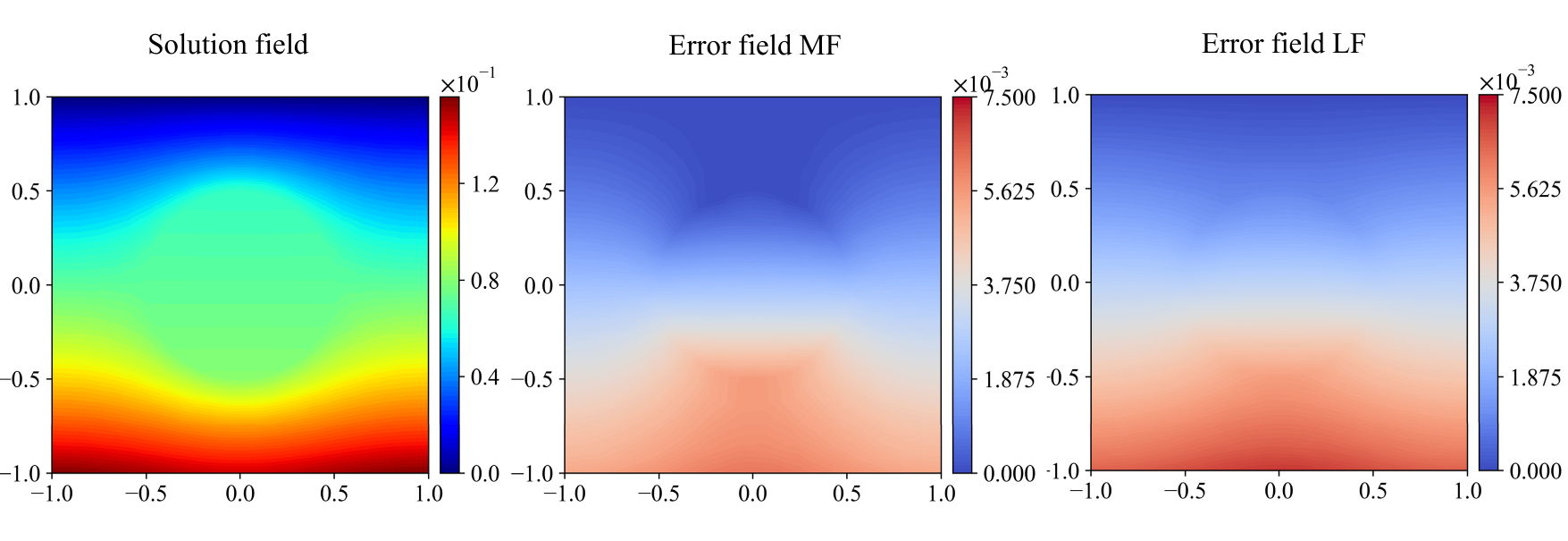}
\end{center}
\vskip-20pt
\caption{Heat conduction problem. From left to right, we depict the high-fidelity solution, and the MF and LF error fields at \((\mu_1, \mu_2, t) = (1.2, -0.1, 3)\).}
\label{fig:fieldsheat}
\end{figure}

Figure \ref{fig:mreheat2} illustrates the timed-averaged relative error ${\bar{\varepsilon}}_{\square}(\bm{\mu})$ distribution in the parametric space. It can be observed that MF-TTOI-DeepONet and LF-TTOI have globally a similar behavior, however, the LF-TTOI exhibits a larger maximum error in the boundary of the parametric space. This further confirms that enriching the information by augmenting the data with surrogates allows DeepONet to learn better the evolution for all parametric instances.  Finally, to evaluate the effectiveness of the proposed methods, Figure \ref{fig:fieldsheat} presents the reference solution and the corresponding error fields of MF-TTOI-DeepONet and LF-TTOI in the spatial domain at final time $t = 3$ for the snapshot with the largest relative error corresponding to $(\mu_1, \mu_2) = (1.2, -0.1)$. The error magnitude is of order $10^{-3}$, confirming the ability of both strategies to achieve accurate generalization, even for parameters outside the training set and in the extrapolation regime.

\subsection{Advection-diffusion problem in a square domain with a moving obstacle}
\label{adeq}
With this second benchmark, we aim at considering a benchmark which actually requires the graph-based architecture due to the more complex and parameterized geometry. Let us consider an advection-diffusion problem in a two-dimensional square domain with a parameterized moving square obstacle \(C(\bm{\mu})\) with a side length of 0.3, so that the computational domain is given by \(\Omega(\bm{\mu}) = [0,1]^2 \setminus C(\bm{\mu})\), with $C(\bm{\mu}) = [\mu_0, \mu_0 + 0.3]\times[\mu_1, \mu_1 + 0.3]$, depicted in Figure \ref{advdomain}.

\definecolor{bluemathlab}{HTML}{065895}
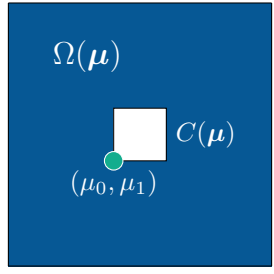
\begin{figure}[htb]
\begin{center}
\centering
\begin{tikzpicture}[scale=3.5]
\fill [bluemathlab] (-0.5,-0.5) rectangle (0.5,0.5);
\fill [white] (-0.1,-0.1) rectangle (0.1,0.1);
\draw[-] (-0.1,-0.1) rectangle (0.1,0.1);
\draw (-0.5,-0.5) -- (0.5,-0.5);
\draw (0.5,-0.5) -- (0.5,0.5);
\draw (0.5,0.5) -- (-0.5,0.5);
\draw (-0.5,0.5) -- (-0.5,-0.5);
\filldraw[fill=greenmathlab, draw=white] (-.1,-.1) circle (1pt);
\node[white] at (-0.2,0.3) {\Large$\Omega(\bm{\mu})$};
\node[white, right] at (0.1,0) {$C(\bm{\mu})$};
\node[white,below] at (-0.1,-0.1) {$(\mu_0,\mu_1)$};

\end{tikzpicture}
\end{center}
\caption{Geometrical setup of the advection-diffusion problem.}
\label{advdomain}
\end{figure}

The physical phenomenon is described by the following governing equation with proper boundary conditions:

$$
\begin{cases}
    \frac{\partial u}{\partial t} - D \Delta u + \mathbf{b} \cdot \nabla u = 0 \quad &\text{in } \Omega(\bm{\mu}) \times (0, t_f],\\
    u(x,y) = (x-1)^2 + (y-1)^2 & \text{on } \partial \Omega(\bm{\mu}) \times (0,t_f], \\
    u(x,y,0) = (x-1)^2 + (y-1)^2 &  \text{in } \Omega(\bm{\mu}),
\end{cases}
$$
where we set \(t_f = 2\), the diffusion coefficient \(D = 0.1\), the advection vector field \(\mathbf{b}(t) = [1 - t, 1 - t]^\top\), and we parameterized the position of the obstacle's center as \(\bm{\mu} = [c_x, c_y] \in \mathcal{P} = [0.1, 0.6]^2\).

The high-fidelity model is obtained using the FE method with \(\mathbb{P}_1\) elements for spatial discretization and the backward Euler method for time integration. The simulations are performed with a time step \(\Delta t = 0.02\), resulting in $N_s=101$ temporal snapshots per simulation.  
For this study, we collected a dataset comprising $N_{\bm{\mu}} = 100$ simulations, considering different obstacle positions, with a total of $N_h = 3159$ mesh nodes is collected. The dataset is split into a training set of $N_\text{train}= 30$ simulations and a testing set of $N_{\text{test}} = 70$ simulations.

\begin{figure}[htbp]
\begin{center}
    \includegraphics[width=0.99\textwidth]{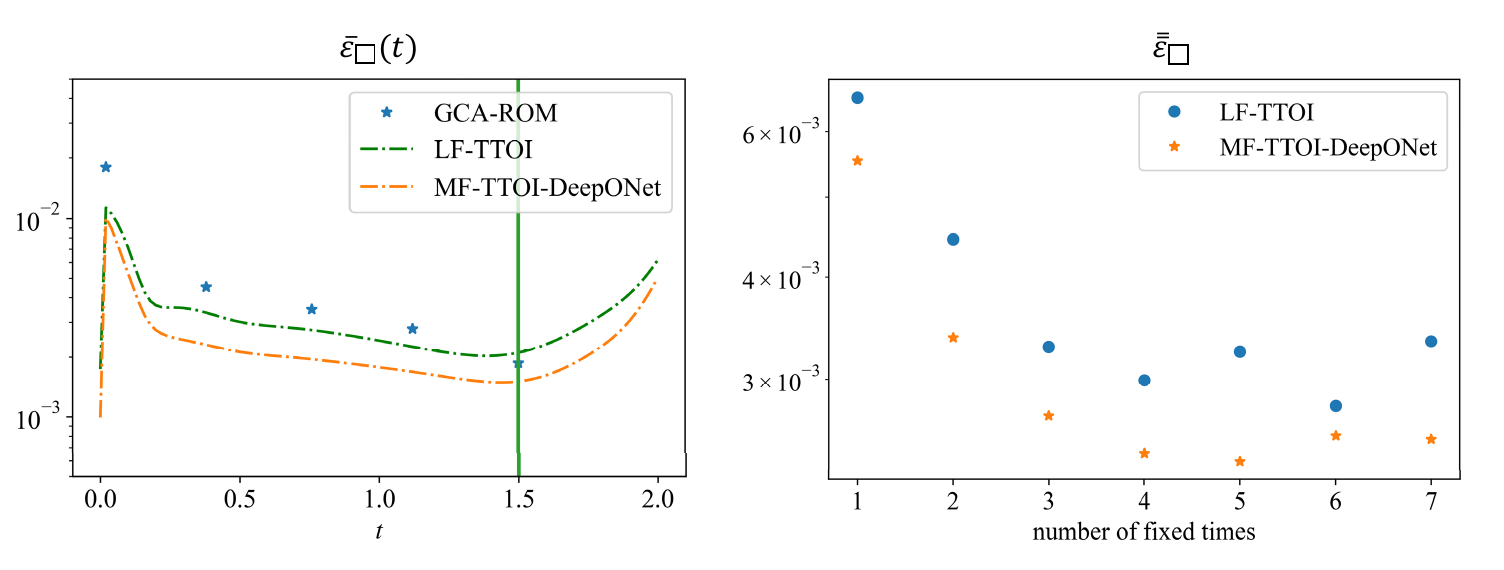}
\end{center}
\vskip-20pt
\caption{Advection-diffusion problem. Left: Mean relative error  ${\bar{\varepsilon}}_{\square}(t)$ of different methods on test parameters over the entire time domain. The vertical line at \( t = 1.5 \) denotes the extrapolation regime. Right: Mean relative error  $\bar{\bar{\varepsilon}}_{\square}$ when exploiting a different number of GCA-ROMs $N_{\text{gca}}$.}
\label{fig:mreadv1}
\end{figure}

A comparison between GCA-ROM, LF-TTOI  and MF-TTOI-DeepONet approaches with $N_{\text{gca}} = 5$ in terms of mean relative error ${\bar{\varepsilon}}_{\square}(t)$ is illustrated in the left plot of Figure \ref{fig:mreadv1}, where the reduced ranks are given by $r_1=30$ and $r_2 = 21$. We can observe that MF-TTOI-DeepONet achieves the best precision on average among the methods, relying on the additional surrogate information to outperform other strategies in the whole time domain. Indeed, incorporating the branch net and learning the residual between the LF-TTOI and the true solution, results in enhanced accuracy.

By varying the amount of GCA-ROM trainings for fixed time instances, i.e.\ $N_{\text{gca}}$, we study the behavior of the error $\bar{\bar{\varepsilon}}_{\square}$ for MF-TTOI-DeepONet and LF-TTOI in the right panel of Figure \ref{fig:mreadv1}. As before, it can be observed that initially, augmenting the number of GCA-ROMs enhances the precision of both methodologies, as a more refined temporal information is available, and until saturation is reached due to limited precision of GCA-ROM.

\begin{figure}[htb]
\begin{center}
    \includegraphics[width=0.99\textwidth]{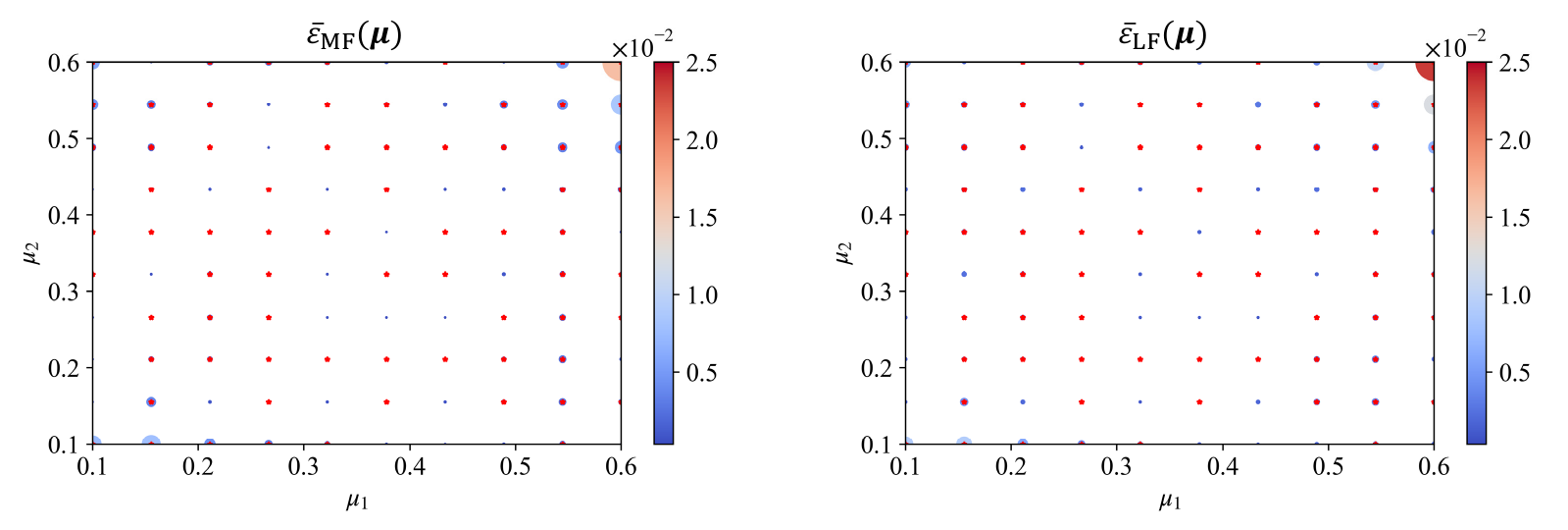}
\end{center}
\vskip-20pt
\caption{Advection-diffusion problem. Mean relative errors  ${\bar{\varepsilon}}_\mathrm{{MF}}(\bm{\mu})$ for MF-TTOI-DeepONet (left) and  ${\bar{\varepsilon}}_\mathrm{{LF}}(\bm{\mu})$ for LF-TTOI (right) on the dataset $\Xi$. Red markers corresponds to the parameters in the testing set $\Xi_{\text{test}}$.}
\label{fig:mreadv2}
\end{figure}

\begin{figure}[htb]
\begin{center}
    \includegraphics[width=0.99\textwidth]{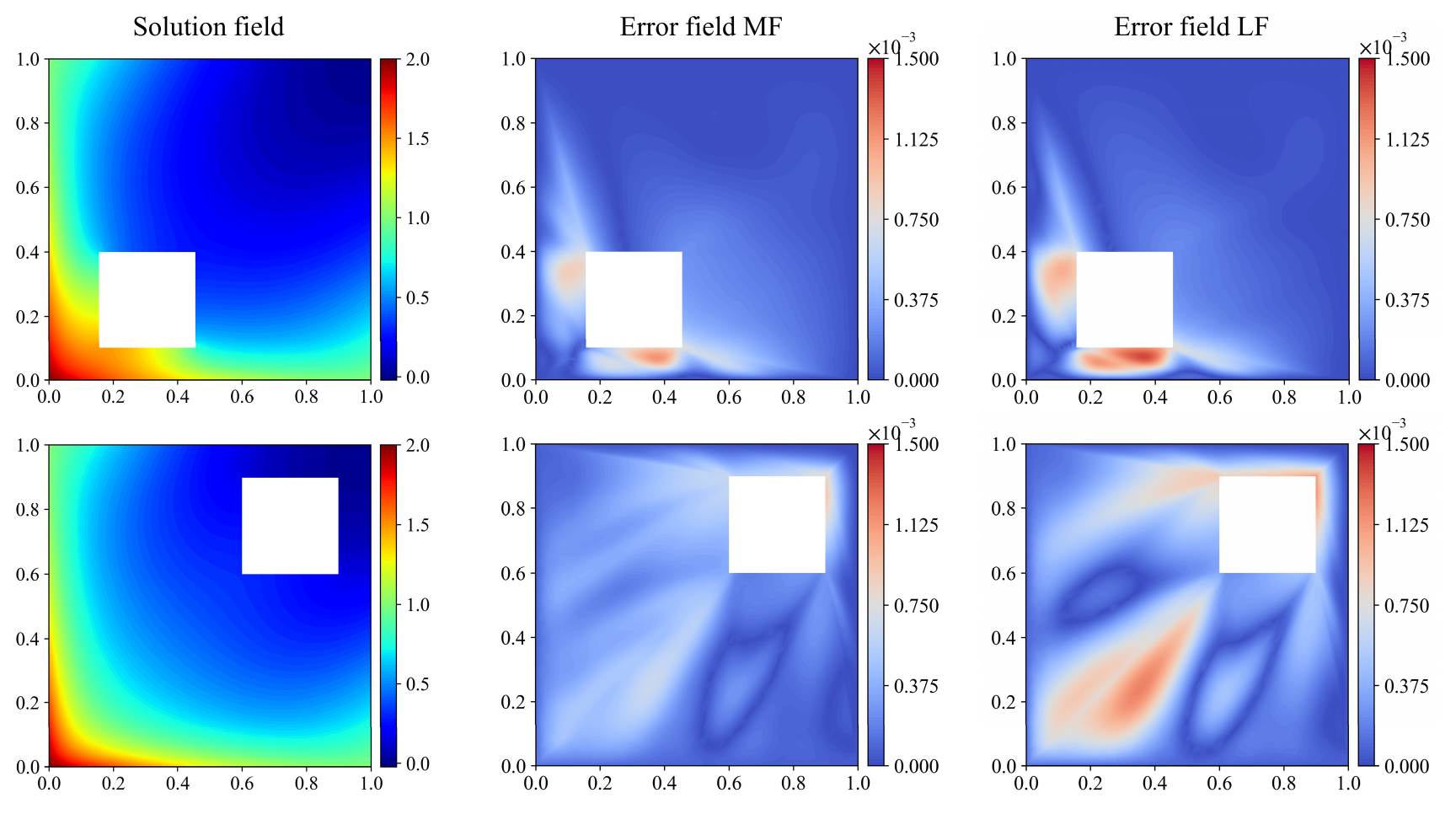}
\end{center}
\vskip-20pt
\caption{Advection-diffusion problem.  The high-fidelity solution, and the MF-TTOI-DeepONet and LF-TTOI error fields. Top: \((\mu_1, \mu_2, t) = (0.1556, 0.1, 2)\).  Bottom: \((\mu_1, \mu_2, t) = (0.6, 0.6, 2)\).}
\label{fig:fieldsadvt}
\end{figure}

To further examine the model's sensitivity to the parameters, Figure \ref{fig:mreadv2} depicts the mean relative errors ${\bar{\varepsilon}}_\mathrm{{MF}}(\bm{\mu})$ and  ${\bar{\varepsilon}}_\mathrm{{LF}}(\bm{\mu})$ for all parameters, additionally illustrating the position of the test samples within the parameter space. The behavior of MF-TTOI-DeepONet and LF-TTOI are similar, except smaller errors are associated with MF-TTOI-DeepONet. This illustrates the advantage of employing the NN architecture to learn the residual. Both distributions exhibit larger errors at the boundary of the parameter domain, which is attributed to the exclusion of these points from the range of the training set, resulting in larger extrapolation errors. Figure \ref{fig:fieldsadvt} illustrates the solution fields and error fields obtained by the MF-TTOI-DeepONet and LF-TTOI at time instance $t = 2$ for the parameter samples $(\mu_1, \mu_2) = (0.1556, 0.10)$ and $(\mu_1, \mu_2) = (0.6, 0.6)$, where both approaches provide an accurate prediction of the solution.

\subsection{Navier-Stokes equations with vortex shedding phenomena}
\label{nseq}
Finally, in this section, we consider the unsteady Navier-Stokes equations in a two-dimensional domain with an obstacle, which describe the evolution of the velocity and pressure fields of an incompressible, viscous fluid. This problem represents an important benchmark in the field of Computational Fluid Dynamics, whose governing equations are given by the following system with proper initial and boundary conditions:

$$
\begin{cases}
    \dfrac{\partial \mathbf{u}}{\partial t} - \nu \nabla^2 \mathbf{u} + (\mathbf{u} \cdot \nabla) \mathbf{u} + \nabla p = 0, & \text{in } \Omega(\bm{\mu}) \times (0, t_f], \\[0.5em]
    \nabla \cdot \mathbf{u} = 0,                                                      & \text{in } \Omega(\bm{\mu}) \times (0, t_f], \\[0.5em]
    \mathbf{u} = 1.5 \left[ \frac{4y(0.41 - y)}{0.41^2}, 0 \right]^\top,                      & \text{on } \Gamma_{\text{in}}\times (0, t_f], \\[0.5em]
    \mathbf{u} = 0,                                                                  & \text{on } \Gamma_w\times (0, t_f], \\[0.5em]
    \mathbf{u} = \mathbf{u}_0,                                                       & \text{in } \Omega(\bm{\mu})\times \{0\}, \\[0.5em]
    \nu \dfrac{\partial \mathbf{u}}{\partial \mathbf{n}} - p\mathbf{n} = 0,          & \text{on } \Gamma_{\text{out}}\times (0, t_f],
\end{cases}
$$
where \(\mathbf{u}\) denotes the velocity field, \( p\) represents the pressure field, and \(\mathbf{n}\) is the unit outward normal vector,  \(\Gamma_{\text{in}}\) is the inflow boundary, \(\Gamma_w\) the no-slip wall boundary, and \(\Gamma_{\text{out}}\) outflow boundary.

\begin{figure}[htb]
\begin{center}
    \includegraphics[width=0.9\textwidth]{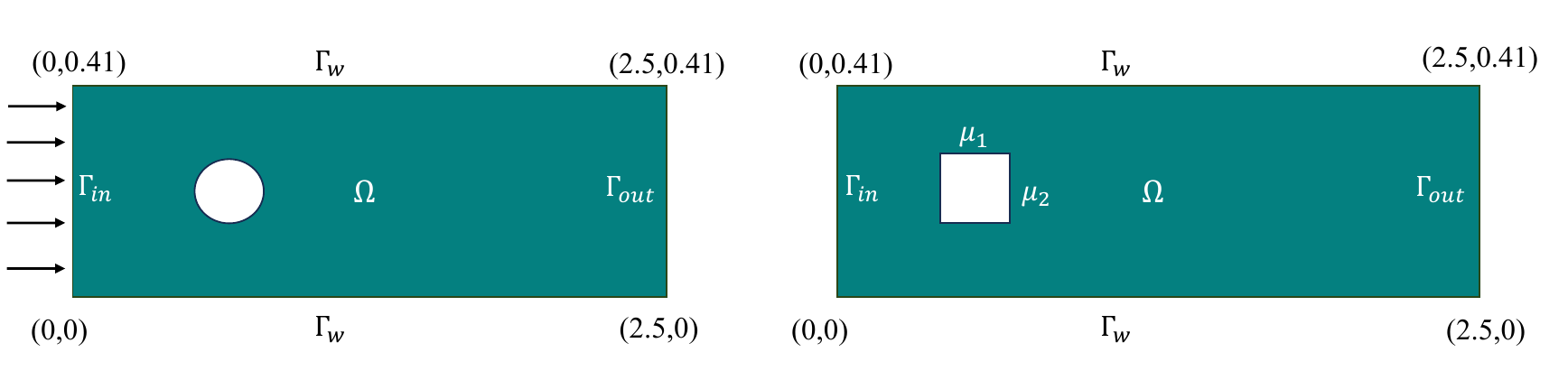}
\end{center}
\vskip-20pt
\caption{Geometrical setup of the unsteady Navier-Stokes problem.}
\label{nsdomain}
\end{figure}

Starting from this formulation, in order to conduct a comprehensive investigation of the phenomenon and evaluate the approximation properties of the proposed strategies, we consider two configurations of the parameterized problem, respectively focusing on the transition between different flow regimes and the shape of the obstacle.
We report now the main features of the two benchmarks of interest:
\begin{itemize}
    \item \textbf{Case 1 (cylinder):} The obstacle is a circular cylinder, as shown in Figure \ref{nsdomain} (left). Since we want to investigate the behavior of the methodology when considering advection-dominated phenomena, i.e.\ flows characterized by high values of the Reynolds number, we consider, as a single parameter, the kinematic viscosity \(\nu\) varying within the range \([3\times10^{-4}, 3\times10^{-3}]\). This corresponds to a Reynolds number \(Re \in [33, 330]\), which is sufficient to observe the transition between different flow regimes. 
    
    \item \textbf{Case 2 (rectangle):} The obstacle is a rectangle centered at $(0.4, 0.205)$, parameterized by its width $\mu_1 \in [0.02, 0.2]$ and height $\mu_2 \in [0.02, 0.2]$, as illustrated in the right of Figure \ref{nsdomain}. In this setting, the kinematic viscosity is fixed at $\nu = 0.0005$, and we aim at revisiting the flow past cylinder benchmark by including a (simple) geometric parametrization of the polygonal obstacle. This could provide interest insights when investigating the lift and drag properties of some shapes immersed in a fluid, or in the field of Optimal Control Problems, when searching for an optimal shape of the obstacle.
\end{itemize}

\subsubsection{Case 1: cylindrical obstacle}\label{case1}
We utilize the FE method with mixed Taylor-Hood \( \mathbb{P}_2-\mathbb{P}_1 \) elements as the high-fidelity model. 
The discretization of the spatial domain results in \( N_h = 5017\) degrees of freedom, while the time derivative approximated by the backward finite difference with a time step \( \Delta t = 0.01 \) up to a final time \( t_f = 6 \).
The high-fidelity dataset is generated by selecting 20 log-equispaced sample points for $\mu = \nu \in \mathbb{P} = [3\times10^{-4}, 3\times10^{-3}]$, with \( N_{\text{train}} = 16 \) randomly selected parameter values used as the training set. By performing TT decomposition, the TT ranks of are $r_1=15$ and $r_2=448$, respectively for the parametric and time dependency, highlighting the complexity of the dynamics transitioning between different states.
The goal is to predict the magnitude velocity $|\mathbf {u} (\mu, t)|$, extrapolating in time, for different values of the kinematic viscosity \(\mu \in \mathbb{P}\) and in the time domain $t \in [0, t_f]$.

\begin{figure}[htb]
\begin{center}
    \includegraphics[width=0.99\textwidth]{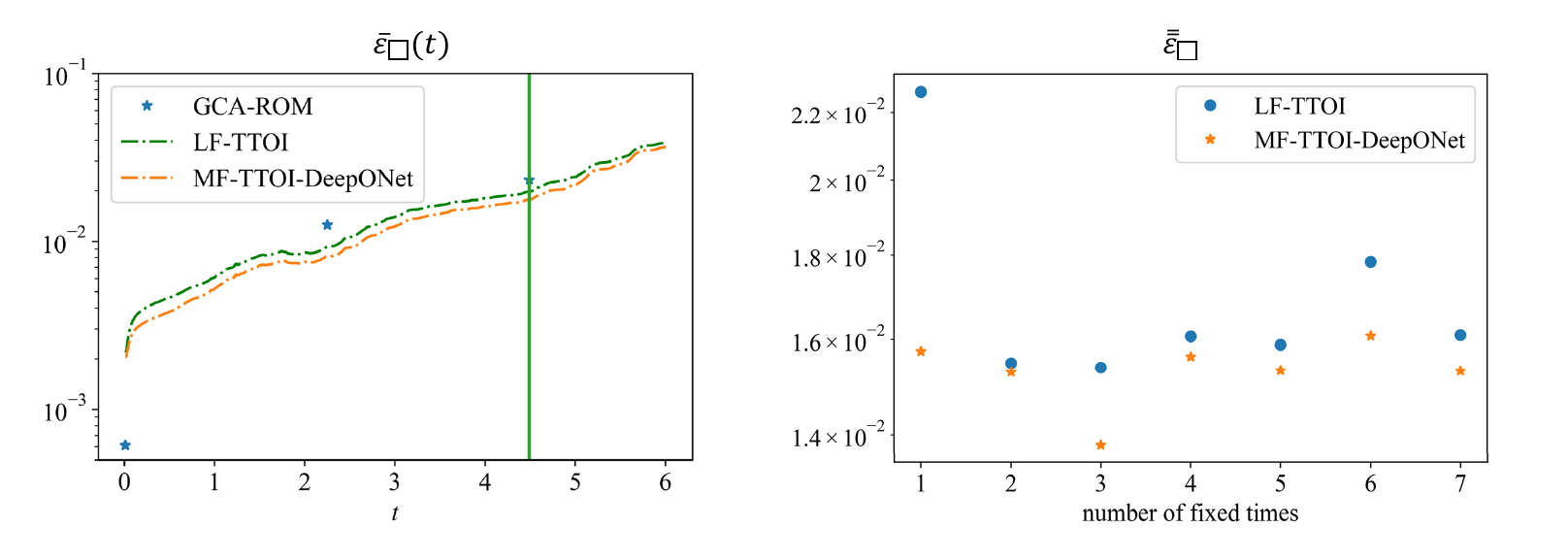}
\end{center}
\vskip-20pt
\caption{Navier-Stokes problem for Case 1. Left:  Mean relative error ${\bar{\varepsilon}}_{\square}(t)$ of different methods over the entire time domain. The vertical line at \( t = 4.49  \) denotes the extrapolation regime. Right: Mean relative error $\bar{\bar{\varepsilon}}_{\square}$ when exploiting a different number of GCA-ROMs $N_{\text{gca}}$.}
\label{fig:mrens1}
\end{figure}

The left plot in Figure \ref{fig:mrens1} presents the mean relative test error ${\bar{\varepsilon}}_{\square}(t)$ for all test parameters over time when the number of GCA-ROMs is fixed at 3.   
MF-TTOI-DeepONet performs a little better than LF-TTOI, especially within the training time range. 
This shows that within the training range time, the branch network can learn the residual between the true solution and the LF-TTOI solution, capable to outperform the other approaches training time range. 
However, the error associated with OpInf becomes the dominant factor, resulting in a similar level of accuracy for the MF-TTOI-DeepONet and LF-TTOI methods.

The trend observed in previous investigations is confirmed by the mean relative error $\bar{\bar{\varepsilon}}_{\square}$ presented in the right plot of Figure \ref{fig:mrens1}, which demonstrates that a relatively small number of GCA-ROMs $N_{\text{gca}}$ is sufficient to provide parameter dependence information, with the selected GCA-ROMs distributed across the entire training time range. We can observe a more pronounced error increase in the extrapolation regime for both strategies w.r.t.\ the previous benchmarks, but we remark this is mostly due to the limited amount of snapshots data computed and the complexity of the phenomena involved.

\begin{figure}[htb]
\begin{center}
    \includegraphics[width=0.99\textwidth]{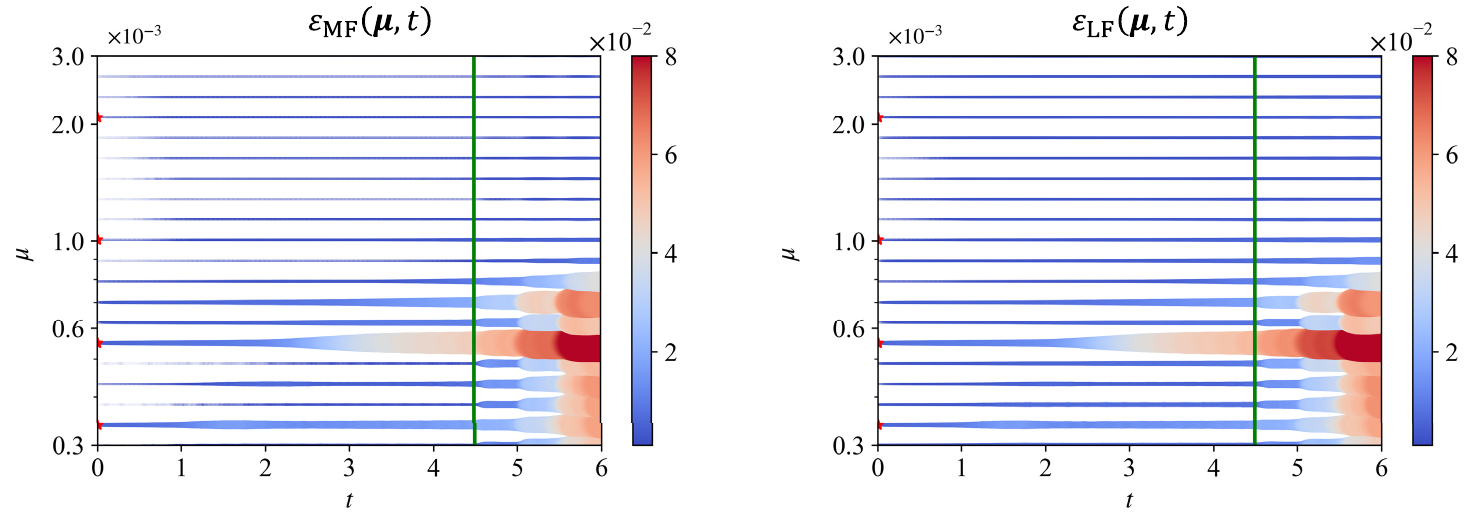}
\end{center}
\vskip-20pt
\caption{Navier-Stokes problem for Case 1. Mean relative error  $\varepsilon_{\text{MF}}(\boldsymbol{\mu},t) $ for MF-TTOI-DeepONet (left) and $\varepsilon_{\text{LF}}(\boldsymbol{\mu},t) $ for LF-TTOI (right). Red markers corresponds to the parameters in the testing set $\Xi_{\text{test}}$.}
\label{fig:mrens2}
\end{figure}

Similarly, we report in Figure \ref{fig:mrens2} the relative errors of MF-TTOI-DeepONet and LF-TTOI for the whole time-domain. 
It can be observed that the errors are concentrated in regions where the viscosity coefficient is relatively low, i.e.\ where convection dominates diffusive phenomena and linear reduction fails to accurately compress the dynamics. The plot also highlights that, with a sufficient amount of sampled data in specific regions of the parametric domain, e.g.\ viscosity values of order $10^{-3}$, the extrapolation predictions in time provide remarkable accuracy even for advection-dominated flows.

\begin{figure}[htb]
\begin{center}
    \includegraphics[width=0.99\textwidth]{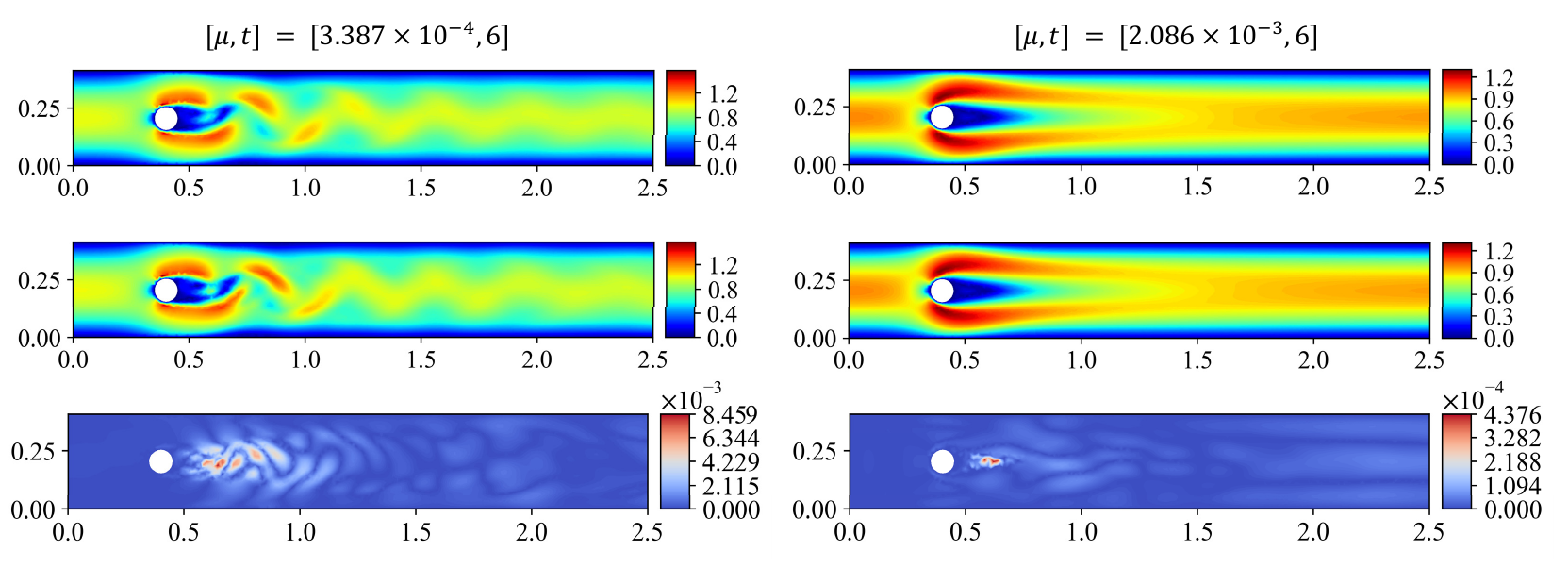}
\end{center}
\vskip-20pt

\caption{Navier-Stokes problem for Case 1.  From top to bottom, we depict the high-fidelity solutions, the MF prediction fields and the error fields at the final time for different parameter samples.}
\label{fig:fieldsns}
\end{figure}

Additionally, we examined the reference solutions, the predicted fields, and their error fields at time $t = 6$ for two different values of the viscosity, \(\mu  = 3.387\times 10^{-4}\) and \(\mu = 2.086\times 10^{-3}\), as illustrated in Figure \ref{fig:fieldsns}. The results demonstrate that the proposed method achieves satisfactory accuracy in capturing the solution characteristics for different flow regimes. 
 
\subsubsection{Case 2: parameterized rectangular obstacle}\label{case2}
For this benchmark, the high-fidelity solutions are computed using the same method as in Section~\ref{case1}, except with a spatial discretization of $N_h = 5424$ and a final simulation time $t_f = 1$. A total of $N_{\bm{\mu}} = 25$ parameter values are sampled through a uniform grid in the parametric space $\mathbb{P} = [0.02,0.2]^2$.  We randomly select $N_\text{train} = 20$ parameter samples for training, while the remaining ones are used for testing, and the training time interval is set to \([0, 0.75]\), the tensor train ranks in this setting are given by $r_1= 20$ and $r_2 = 76$, denoting, as expected, a more complex parametrization and an easier dynamical evolution.

\begin{figure}[htb]
\begin{center}
    \includegraphics[width=0.99\textwidth]{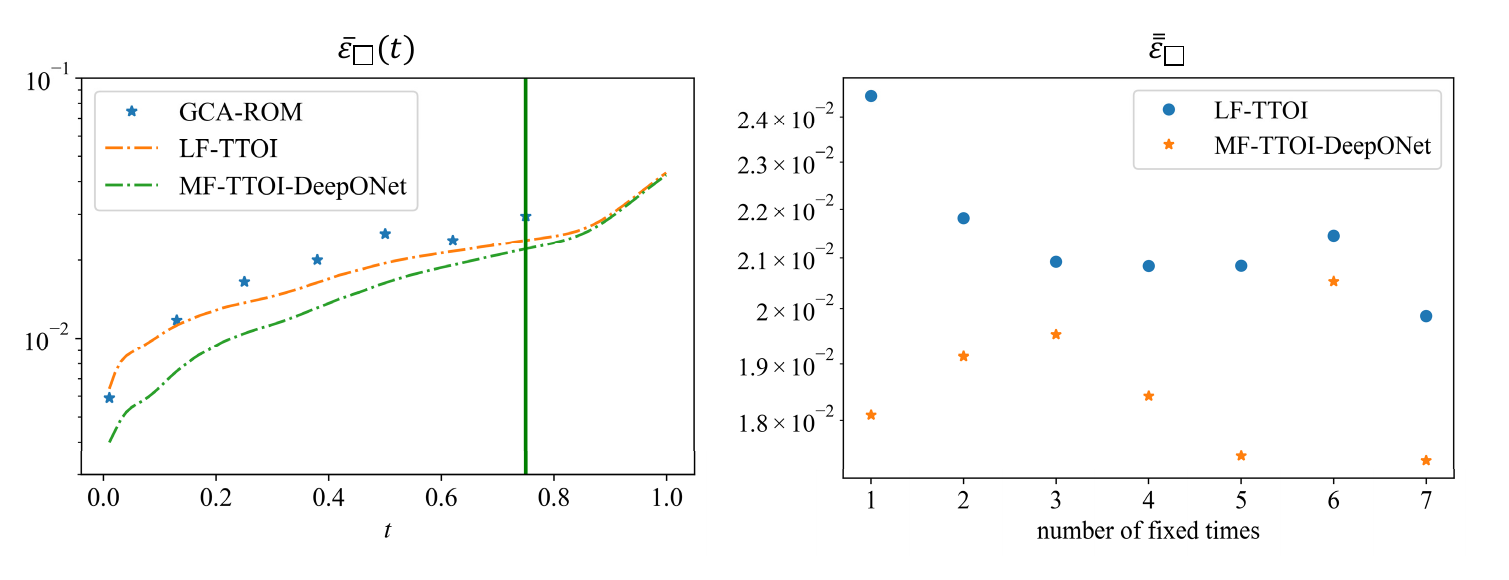}
\end{center}
\vskip-20pt
\caption{Navier-Stokes problem for Case 2. Left: Mean relative error ${\bar{\varepsilon}}_{\square}(t)$ of different methods over the entire time domain.  The vertical line at \( t = 0.75  \) denotes the extrapolation regime. Right: Mean relative error $\bar{\bar{\varepsilon}}_{\square}$ when exploiting a different number of GCA-ROMs $N_{\text{gca}}$.}
\label{fig:mrens2moving}
\end{figure}

By selecting $N_{\text{gca}} = 7$, we show the mean relative errors ${\bar{\varepsilon}}_{\square}(t)$ of  LF-TTOI  and MF-TTOI-DeepONet in the left panel of Figure \ref{fig:mrens2moving}. The results show that MF-TTOI-DeepONet consistently outperforms LF-TTOI  across all time instances, achieving a maximum mean relative error of order $4 \times 10^{-2}$.  Once again, we depict in the right panel of Figure \ref{fig:mrens2moving} the impact of the number of GCA-ROMs on model accuracy, for LF-TTOI  and MF-TTOI-DeepONet, by means of the mean relative errors $\bar{\bar{\varepsilon}}_{\square}$. 

It can be observed that the proposed method achieves a satisfactory approximation already with as few as three GCA-ROMs, stagnating due to the error from GCA-ROM and OpInf time prediction, and from the lack of many additional test parameters to augment the information and benefit from the two-stage strategy.

\begin{figure}[htb]
\begin{center}
    \includegraphics[width=0.99\textwidth]{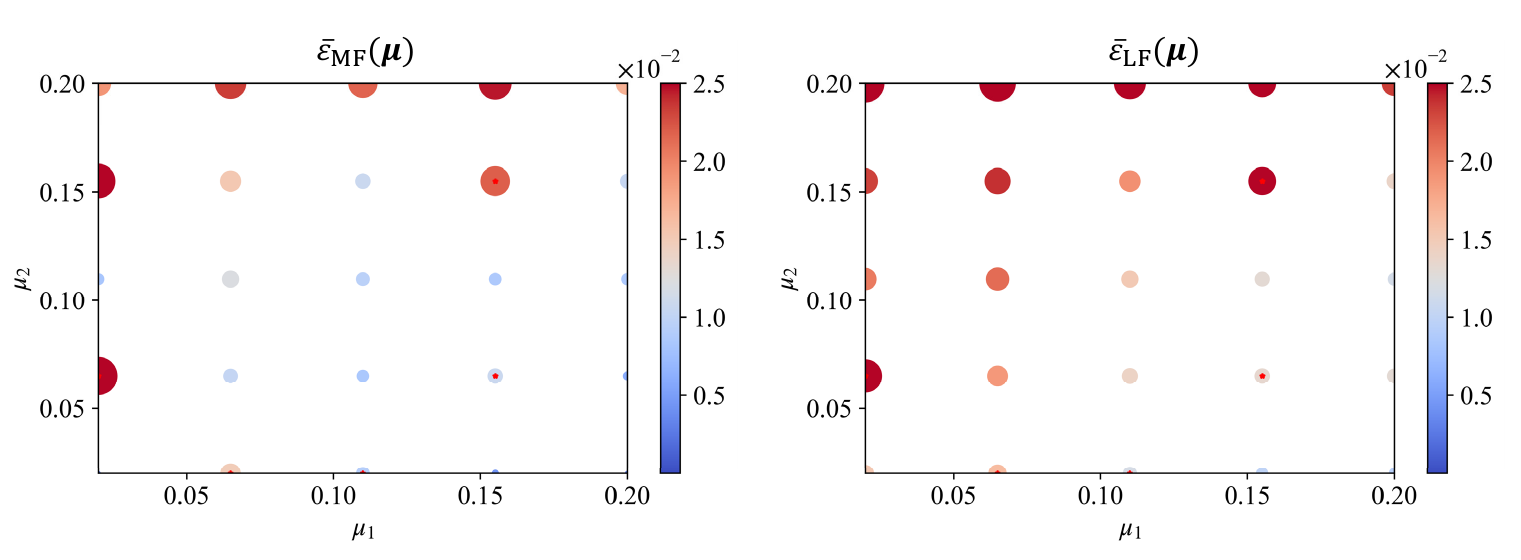}
\end{center}
\vskip-20pt
\caption{Navier-Stokes problem for Case 2. Mean relative errors ${\bar{\varepsilon}}_{\mathrm{MF}}(\bm{\mu})$ for MF-TTOI-DeepONet (left) and ${\bar{\varepsilon}}_{\mathrm{LF}}(t)$ for LF-TTOI (right) on the dataset $\Xi$. Red markers corresponds to the parameters in the testing set $\Xi_{\text{test}}$.}
\label{fig:mrens1moving}
\end{figure}

To examine the predictive performance of the model with respect to different parameter values, we plotted the time-averaged relative errors ${\bar{\varepsilon}}_{{\square}}(\bm{\mu})$ for various parameters in Figure \ref{fig:mrens1moving}. The performance of  LF-TTOI  and MF-TTOI-DeepONet appears largely similar across geometric parameters, with MF-TTOI-DeepONet exhibiting higher accuracy, as expected. The error distribution further indicates that the model is more sensitive to the geometrical parameter \( \mu_2 \), which creates very different flow behaviors between the obstacle and the walls.

\begin{figure}[htb]
\begin{center}
    \includegraphics[width=0.99\textwidth]{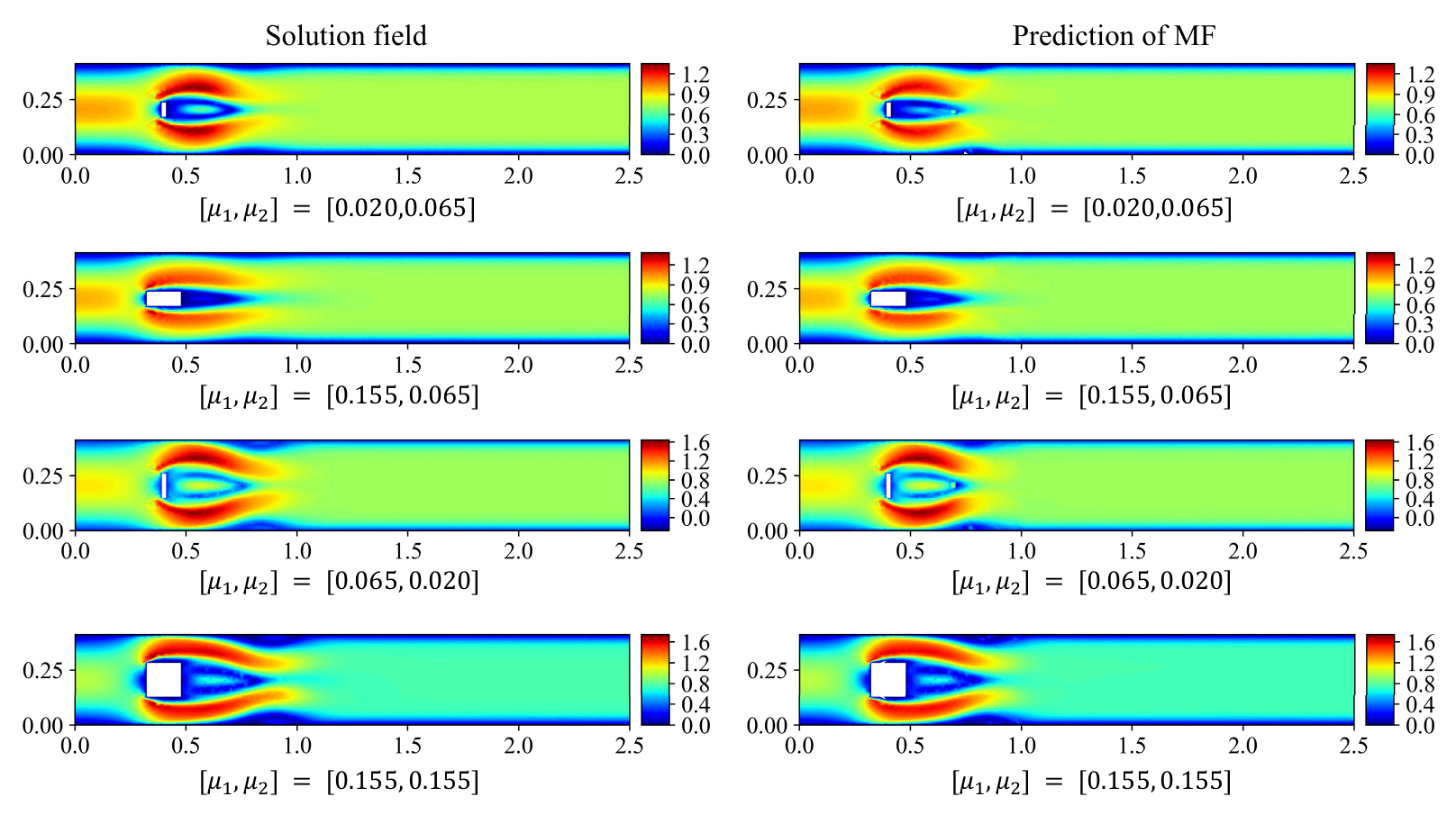}
\end{center}
\vskip-20pt

\caption{Navier-Stokes problem for Case 2. From left to right, we depict the high-fidelity solutions and the MF prediction fields at time $t = 1$ for different parameter samples.}
\label{fig:fieldsnsmoving}
\end{figure}

To provide a more detailed analysis, we selected four distinct parameter sets that comprehensively illustrate how the solution varies with parameter changes, and visualized both the reference solutions at the final time step and the corresponding predictions from the MF-TTOI-DeepONet in Figure \ref{fig:fieldsnsmoving}. The spatial distribution confirms that the proposed model effectively captures the spatial characteristics of the solution across different parameter settings.

\subsection{A comprehensive comparative study on extrapolation and methodologies}
In this section, we evaluate the performance of the proposed methods across different extrapolation regimes, i.e.\ training time ranges, by computing the mean relative errors of each method, as presented in Table \ref{mretime}.  We remark that we do not focus on selecting the optimal GCA-ROM for each case. Instead, we consistently used the same number of GCA-ROM models in all time range, so the performance of the proposed models may not necessarily be optimal. Our primary focus is on evaluating the applicability of the proposed models for longer time predictions, when only a limited amount of training data is available. The results show that in the initial stages, the error of the proposed methods decreases significantly as the training time range increases and then gradually stabilizes. This phenomenon can be explained by the fact that when the training time range is short and the prediction horizon is long, the OpInf model fails to adequately capture the system's temporal evolution, and such error dominates the one of the whole architecture. As the training time range increases, the accuracy of the OpInf model improves, leading to a substantial reduction in the overall error. However, further extending the training time range does not result in additional improvement because the estimation accuracy of the parameter coefficients becomes the primary limiting factor, and parameter estimation errors dominate at this stage.

\begin{table}[h!] 
\centering
\begin{tabular}{cccccc}
\hline
 \textbf{Application} & \textbf{Method}  &  [0, $\frac{1}{4}t_f$]       & [0, $\frac{1}{2}t_f$]       & [0, $\frac{3}{4}t_f$]       & [0, $t_f$] \\ \hline
\multirow{2}{*}{Heat} 
& MF & $1.07 \times 10^{-1}$ & $6.07 \times 10^{-3}$ & $6.26 \times 10^{-3}$ & $6.46 \times 10^{-3}$ \\ 
& LF & $1.14 \times 10^{-1}$ & $7.35 \times 10^{-3}$ & $7.53 \times 10^{-3}$ & $8.37 \times 10^{-3}$ \\ 

\multirow{2}{*}{Advection-diffusion} 
& MF & $1.07 \times 10^{-1}$ & $1.29 \times 10^{-2}$ & $2.38 \times 10^{-3}$ & $2.35 \times 10^{-3}$ \\ 
& LF & $1.08 \times 10^{-1}$ & $1.38 \times 10^{-2}$ & $3.24 \times 10^{-3}$ & $3.24 \times 10^{-3}$ \\ 

\multirow{2}{*}{Navier-Stokes Case1}  
& MF & $2.36 \times 10^{0}$ & $7.59 \times 10^{-2}$ & $1.46 \times 10^{-2}$ & $1.60 \times 10^{-2}$ \\ 
& LF & $2.37 \times 10^{0}$ & $7.68 \times 10^{-2}$ & $1.59 \times 10^{-2}$ & $1.66 \times 10^{-2}$ \\ 

\multirow{2}{*}{Navier-Stokes Case2}  
& MF & $6.48 \times 10^{-2}$ & $2.97 \times 10^{-2}$ & $1.73 \times 10^{-2}$ & $1.67 \times 10^{-2}$ \\ 
& LF & $6.52 \times 10^{-2}$ & $3.06 \times 10^{-2}$ & $1.99 \times 10^{-2}$ & $1.86 \times 10^{-2}$ \\ \hline
\end{tabular}
\caption{The mean relative errors of different training time ranges for the MF-TTOI-DeepONet and LF-TTOI.}
\label{mretime}
\end{table}

Moreover, we compare the proposed methods with non-intrusive state-of-the-art reduced-order approaches: (i) POD-Proj, projecting high-dimensional snapshots onto a POD basis subspace as the optimal linear approximation; (ii) POD-NN \cite{HESTHAVEN201855}, which employs neural networks to map parameters to projection coefficients for nonlinear dynamics; and (iii) MeshGraphNets \cite{pfaff2020learning}, that operates directly on parameterized mesh graphs, using both geometry and current state as input to predict future states.
By setting the energy threshold for the POD to $10^{-8}$, consistently with the TT decomposition of the previous analysis, the number of reduced-basis of the POD-based methods for Examples \ref{heateq}, \ref{adeq} and the two cases of \ref{nseq} are 11, 82, 350 and 295 respectively.

As concerns MeshGraphNets, since the methodology consists of an autoregressive approach, it is particularly well suited for predicting dynamics defined on complex unstructured grids. Despite this, the original version is not particularly efficient to treat physically parameterized problems, since the input only includes the current state along with mesh information but no parameter vectors. 

For this reason, we compare MeshGraphNets only to the ``purely geometrical" benchmarks presented above, the advection-diffusion problem in Section \ref{adeq} and the Navier-Stokes problems with rectangular obstacle in Section \ref{case2}. Table \ref{mremethod} the accuracy of the different methodologies when using the same set of snapshots. The table shows that, as expected, when a sufficiently large number of bases is selected the best linear projection error (not useful for online prediction) outperforms all strategies, but a simple regression task based on such basis to recover the reduced coefficients via Neural Network do not provide accurate results. 

For all benchmark problems, the errors of MF-TTOI-DeepONet and LF-TTOI are significantly lower than those of the POD-NN method, with MF-TTOI-DeepONet demonstrating an improvement of the accuracy. For the simple heat conduction problem, the POD-NN method struggled to produce accurate estimates due to the insufficient temporal training data. However, both LF-TTOI and MF-TTOI-DeepONet achieved high accuracy, benefiting from the OpInf method's ability to effectively capture the system's temporal evolution. Furthermore, for more complex problems, the proposed methods maintained superior performance, attributed to the excellent performance of GCA-ROM with limited training parameters and the extrapolation capability of the OpInf approach.

From the comparison with MeshGraphNets models, we can observe that, for both advection-diffusion and Navier-Stokes problems, the existing training data was not sufficiently big to train an accurate surrogate model, demonstrating the potentiality of our strategy to be data-efficient while capable of accurately predicting extrapolation results.

\begin{table}[h!]
\centering
\resizebox{\textwidth}{!}{%
\begin{tabular}{@{}lccccc@{}}
\toprule
\textbf{Application}       &  \textbf{POD-Proj} &  \textbf{POD-NN}      & \textbf{LF-TTOI}      & \textbf{MF-TTOI-DeepONet}      & \textbf{MeshGraphNets}   \\ \midrule
Heat                    & $1.50 \times 10^{-4}$ & $2.22 \times 10^{-1}$ & $8.78 \times 10^{-3}$ & $7.63 \times 10^{-3}$ & -- \\
Advection-diffusion     & $5.69 \times 10^{-3}$ & $8.33 \times 10^{-2}$ & $4.15 \times 10^{-3}$ & $2.51 \times 10^{-3}$ & $4.09 \times 10^{-1}$ \\
Navier-Stokes Case 1    & $4.50 \times 10^{-4}$ & $5.28 \times 10^{-2}$ & $1.54 \times 10^{-2}$ & $1.38 \times 10^{-2}$ & -- \\
Navier-Stokes Case 2    & $2.47 \times 10^{-3}$ & $8.64 \times 10^{-2}$ & $1.99 \times 10^{-2}$ & $1.73 \times 10^{-2}$ & $2.12 \times 10^{-1}$ \\ \bottomrule
\end{tabular}
}
\caption{The mean relative errors for different techniques on the benchmark problems.}
\label{mremethod}
\end{table}

\section{Conclusion}\label{Conclusion}
This work first introduces a novel data-driven reduced-order modeling approach for parameterized time-dependent PDEs, LF-TTOI, integrating GCA-ROM and Operator Inference for parameter and time reduction in a comprehensive framework enabled by Tensor Train decomposition. In the offline stage, the high-fidelity snapshots are approximated as tensor products of parametric, spatial, and temporal components by performing the TT decomposition on the collected snapshot tensor. OpInf learns the temporal behavior of the system from the time-dependent data in order to provide long-term predictions on the future evolution of the time cores. Several parameterized GCA-ROMs at fixed times are trained to provide parameter-dependent information on the cores. In the online stage, OpInf is used to estimate the cores at the new time, while the cores at the new parameters can be computed by minimizing the error between the prediction of the GCA-ROMs at the new parameters and the low-rank approximation. In addition, we developed a multi-fidelity modeling framework, MF-TTOI-DeepONet, in which the combination of temporal and spatial cores corresponds to the trunk net. The branch net is used to predict the parameter cores, and the residual learning is performed by adding the coefficients calculated by LF-TTOI to the input to further improve the prediction accuracy. The predictions of GCA-ROMs and collected high-fidelity snapshots are used to co-train the branch net.

The proposed methods are tested for complex physical scenarios with parameterized geometries to efficiently extrapolate in-time, from heat conduction to Navier-Stokes problem via advection-diffusion benchmark, where both physical (e.g.\ viscosity) and geometrical (e.g.\ obstacles dimension/position) parameters are treated. The numerical results demonstrate the ability of the proposed strategies to achieve higher accuracy than the (causality-lacking) GCA-ROM approach only exploiting a limited number of temporal slices, and also exhibits excellent extrapolation capability for systems with varying parameters. Furthermore, the multi-fidelity approach outperforms the LF-TTOI  by leveraging a DeepONet structure to augment data in a multi-fidelity paradigm, enhance parameter coefficient estimations, and achieving remarkable accuracy even in the low-data regime. Additionally, the frameworks introduced in this study are not confined to GCA-ROM and can be extended to a broader range of data-driven models.

In future work, several research directions can be considered. First, we aim to investigate benchmarks for even challenging cases, including large scale problems and high-dimensional parametric spaces. To further reduce the computational cost of the DeepONet, adaptive sampling strategies could be explored, optimizing the loss function using only a subset of spatial points. Second, the integration of different discretization techniques with operator modeling will be studied to improve the generality of the proposed methods. Interesting directions could also involve the exploitation of the same framework but adopting novel strategies both for the parametric dependence such as Graph Feedforward Network (GFN) \cite{MORRISON2024117458}, and the time behavior via Neural ODEs \cite{ChenNeuralOrdinaryDifferential2018} and Latent Dynamics Networks \cite{RegazzoniLearningIntrinsicDynamics2024}.

\section*{Acknowledgements}
\noindent \textbf{FP} and \textbf{GR} acknowledge the support provided by the European Union - NextGenerationEU, in the framework of the iNEST - Interconnected Nord-Est Innovation Ecosystem (iNEST ECS00000043 - CUP: G93C22000610007) consortium and its CC5 Young Researchers initiative.  \textbf{ZG} and \textbf{YC}  acknowledge the support provided by the China Scholarship Council (No. 202406330108), the National Natural Science Foundation of China (12371435), the Taishan Scholars Program (tsqn202211059) and the Shandong Provincial Natural Science Foundation (ZR2023MA043). The authors also like to acknowledge INdAM-GNSC and the project ``Sviluppo e analisi di modelli di ordine ridotto basati su tecniche di deep learning" (CUP: E53C24001950001) for its support.

\bibliographystyle{elsarticle-num}
\bibliography{main}%

\end{document}